\documentclass{newbaki}
\usepackage{xspace,smfenum}
\usepackage{mathrsfs,amssymb}
\usepackage[T1]{fontenc}
\usepackage[frenchb]{babel}
\usepackage[cm]{aeguill}

\date{Mars 2001}
\bbkannee{53\textsuperscript{e} ann{\'e}e, 2000-2001} 
\bbknumero{886}

\def\A{{\mathbf A}}
\def\Z{{\mathbf Z}}
\def\R{{\mathbf R}}
\def\Q{{\mathbf Q}}
\def\N{{\mathbf N}}
\def\C{{\mathbf C}}
\def\F{{\mathbf F}}
\def\P{{\mathbf P}}
\def\Gm{{{\mathbf G}_m}}

\def\Jac{\operatorname{Jac}}
\let\emptyset\varnothing
\def\Lie{\operatorname{Lie}}
\def\Spec{\operatorname{Spec}}
\def\hdeg{\mathop{\widehat{\mathrm {deg}}}}
\def\hmu{\mathop{\widehat{\mu}}}
\def\hmumax{\mathop{\widehat{\mu}_{\max}}}
\def\rg{\operatorname{rang}}
\let\ra\rightarrow
\let\hra\hookrightarrow
\def\abs#1{\left\lvert {#1} \right\rvert}
\def\norm#1{\left\lVert {#1} \right\rVert}
\def\transp #1{\vphantom{#1}^{\mathrm t}\! {#1}}
\def\CEnd{\operatorname{\mathit{End}}}
\def\End{\operatorname{End}}
\def\GL{\operatorname{GL}}

\def\Sym{\operatorname{Sym}}
 
\def\bar#1{\overline{#1}\vphantom{#1}}
\def\pieme{$p$-i{\`e}me\xspace}

\let\leq\leqslant
\let\geq\geqslant
\let\phi\varphi
\let\eps\varepsilon

\makeatletter
\def\varlimsup@#1{\mathop{\@@overline
  {\hbox{$#1\m@th\operator@font lim\vphantom{fp}$}}}}
\makeatother
\let\limsup\varlimsup

\let\mathcal\mathscr

\def\labelenumi{\arabic{enumi}\textsuperscript o)}

\swapnumbers

\theoremstyle{plain}
\newtheorem{Theo}[subsection]{Th{\'e}or{\`e}me}
\newtheorem{Conj}[subsection]{Conjecture}
\newtheorem{Prop}[subsection]{Proposition}
\newtheorem{Lemm}[subsection]{Lemme}

\theoremstyle{remark}
\newtheorem{Rema}[subsection]{Remarque}

 \makeatletter

 \def\subsection{\@startsection{subsection}{2}%
   \z@{.7\linespacing\@plus.3\linespacing}{-0pt}
   {\normalfont\bfseries\smf@boldmath}}

 \let\c@equation\c@subsubsection
 \let\cl@equation\cl@subsubsection
 \let\theequation\thesubsubsection
 \makeatother
 
\begin{document}
\frontmatter

\title[Th{\'e}or{\`e}mes d'alg{\'e}bricit{\'e}]
      {Th{\'e}or{\`e}mes d'alg{\'e}bricit{\'e}\\ en g{\'e}om{\'e}trie diophantienne}
\subtitle{d'apr{\`e}s J.-B.~Bost, Y.~Andr{\'e}, D.\ \&\ G.~Chudnovsky}
\alttitle{Algebraicity theorems in diophantine geometry}

\author[A. Chambert-Loir]{Antoine CHAMBERT-LOIR}
\address{Centre de Math{\'e}matiques\\
{\'E}cole polytechnique \\
F-91128 Palaiseau Cedex}
\email{chambert@math.polytechnique.fr}

\begin{abstract}
Dans un article r\'ecent, J.-B.~Bost 
\'etablit un crit\`ere
assurant que certaines {\og sous-vari\'et\'es formelles\fg}
de vari\'et\'es alg\'ebriques sont en fait alg\'ebriques.
Son th\'eor\`eme unifie et g\'en\'eralise des r\'esultats des
fr\`eres Chudnovsky et de Y.~Andr\'e motiv\'es par une conjecture
arithm\'etique de Grothendieck pr\'edisant que les solutions
de certaines \'equations diff\'erentielles sont des fonctions
alg\'ebriques.

La d\'emonstration reprend les techniques d'approximation diophantienne
utilis\'ees par ces auteurs avec toutefois un point de vue
syst\'ematiquement g\'eom\'etrique, via notamment
la g\'eom\'etrie d'Arakelov et le formalisme des pentes.
\end{abstract}

\begin{altabstract}
In a recent paper, J.-B. Bost proves a criterion which insures
that certain ``formal subvarieties'' of algebraic varieties are
in fact algebraic. His theorem unifies and generalizes previous
results by D. and G. Chudnovksy, and by Y.~Andr\'e.
They were motivated
by an arithmetic conjecture of Grothendieck which predicts
that the solutions of certain differential equations are
algebraic functions.

The proof uses diophantine approximation techniques
pioneered in this context by these authors, but with
a systematic geometrical point of view, notably using
Arakelov geometry and the formalism of slopes.
\end{altabstract}


\altkeywords{Arakelov geometry, slope method, algebraic leaves of foliations,
  arithmetic differential equations, Liouville property} 

\maketitle
\mainmatter

\section{Introduction}\label{sec.intro}

\subsection{}
Certains probl{\`e}mes g{\'e}om{\'e}triques ou arithm{\'e}tiques
se ram{\`e}nent parfois {\`a} d{\'e}cider si une fonction
analytique, voire une s{\'e}rie formelle, est en fait une fonction
alg{\'e}brique, voire m{\^e}me une fraction rationnelle.
On peut notamment penser {\`a} deux exemples :
\begin{enumerate}
\item le probl{\`e}me de Schwarz (cf.~\cite{schwarz1873})
consistant {\`a} d{\'e}terminer les {\'e}quations
diff{\'e}rentielles hyperg{\'e}om{\'e}triques de Gau{\ss} ayant une base de solutions
form{\'e}e de fonctions alg{\'e}briques ;
\item le th{\'e}or{\`e}me de {\'E}.~Borel~\cite{eborel1894} affirmant qu'une s{\'e}rie enti{\`e}re
$f\in \Z[[x]]$ est le d{\'e}veloppement en s{\'e}rie d'une fraction
rationnelle si et seulement si c'est le d{\'e}veloppement
de Taylor en l'origine d'une fonction \emph{m{\'e}romorphe}
sur un  disque de centre $0$ et rayon strictement sup{\'e}rieur {\`a} $1$.
Une g{\'e}n{\'e}ralisation de ce crit{\`e}re a {\'e}t{\'e} utilis{\'e}e par Dwork~\cite{dwork60}
dans sa d{\'e}monstration de la rationalit{\'e} de la fonction z{\^e}ta d'une
vari{\'e}t{\'e} alg{\'e}brique sur un corps fini.
\end{enumerate}

Les travaux de D.V. et G.V.~Chudnovsky~\cite{chudnovsky85b}
ont montr{\'e} que les m{\'e}thodes de transcendance fournissent
une approche {\`a} ce genre de questions, en m{\^e}me temps qu'{\`a}
d'autres probl{\`e}mes a priori plus {\'e}loign{\'e}s telle
la conjecture d'isog{\'e}nie. C'est l{\`a} le sujet de notre expos{\'e}.
Nous en donnerons de nombreuses applications, dont
certaines remontent d'ailleurs {\`a} ces auteurs.
Cette technique a toutefois {\'e}t{\'e} perfectionn{\'e}e, notamment
par Y.~Andr{\'e} (\cite{andre89} et~\cite{andre97}) et,
plus r{\'e}cemment, par un article de J.-B. Bost~\cite{bost2001}
qui fournit un crit{\`e}re
permettant d'affirmer que certaines sous-vari{\'e}t{\'e}s formelles
d'une vari{\'e}t{\'e} alg{\'e}brique 
(c'est-{\`a}-dire {\og d{\'e}finies\fg} par des s{\'e}ries formelles plut{\^o}t
que des polyn{\^o}mes) sont en fait des sous-vari{\'e}t{\'e}s alg{\'e}briques.
Celui-ci est {\'e}crit dans le langage de la g{\'e}om{\'e}trie d'Arakelov
et j'esp{\`e}re convaincre le lecteur qu'elle fournit un cadre
g{\'e}om{\'e}trique efficace pour les d{\'e}monstrations d'approximation
diophantienne.

Les th{\'e}or{\`e}mes ant{\'e}rieurs
des  Chudnovsky et Andr{\'e} (voir
% de D.~V. et G. V.~Chudnovsky et Y.~Andr{\'e} (voir
notamment~\cite{chudnovsky85b,andre89,andre97})
{\'e}tablissaient l'alg{\'e}bricit{\'e} d'une s{\'e}rie formelle
et les cons{\'e}quences g{\'e}om{\'e}triques {\'e}taient alors obtenues apr{\`e}s
d{\'e}vissage. Les r{\'e}sultats de Bost g{\'e}n{\'e}ralisent 
aussi un th{\'e}or{\`e}me 
de Graftieaux (voir~\cite{graftieaux2001,graftieaux2001b})
concernant les sous-groupes formels d'une vari{\'e}t{\'e} ab{\'e}lienne,
tout en en simplifiant notablement la d{\'e}monstration.

% \subsection{}
% Quel que soit le langage utilis{\'e} (approximation diophantienne
% classique, g{\'e}om{\'e}trie d'Arakelov,\dots), la {\og m{\'e}thode de Chudnovsky\fg}
% repose sur la conjonction de deux types d'hypoth{\`e}ses:
% \begin{itemize}
% \item pour presque tout nombre premier $p$,
% on suppose que 
% la situation est $p$-adiquement bien meilleure que ce qu'on
% pourrait pr{\'e}voir. Ce sont ainsi des hypoth{\`e}ses d'int{\'e}gralit{\'e}, 
% ou plus g{\'e}n{\'e}ralement un contr{\^o}le des d{\'e}nominateurs ;
% \item en une place complexe, on se donne une uniformisation
% holomorphe de la sous-vari{\'e}t{\'e} formelle par une vari{\'e}t{\'e}
% complexe {\og assez grosse\fg}. Dans les articles des Chudnovsky,
% c'est par exemple une uniformisation par une fonction
% sur un espace affine $\C^N$ d'ordre exponentiel fini. Dans l'article de Bost,
% cette condition est affaiblie et il requiert l'existence
% d'une \emph{uniformisation} par une vari{\'e}t{\'e} analytique
% complexe~$M$ telle que l'espace affine,
% ou, plus g{\'e}n{\'e}ralement, satisfaisant une condition dite  \emph{de Liouville},
% {\`a} savoir que toute fonction plurisousharmonique born{\'e}e sur~$M$
% est constante. 
% \end{itemize}
% 
% Pour tirer des cons{\'e}quences arithm{\'e}tiques de cette m{\'e}thode,
% il est en g{\'e}n{\'e}ral n{\'e}cessaire de se d{\'e}faire de cette derni{\`e}re
% hypoth{\`e}se. Dans certaines situations, notamment celles o{\`u}
% interviennent des groupes alg{\'e}briques commutatifs, l'application
% exponentielle fournit l'uniformisation voulue.

\subsection{}
Quel que soit le langage utilis{\'e} (approximation diophantienne
classique, g{\'e}om{\'e}trie d'Arakelov,\dots), la {\og m{\'e}thode de Chudnovsky\fg}
repose sur des techniques utilis{\'e}es en th{\'e}orie des nombres
transcendants ; elle combine en effet
des hypoth{\`e}ses de type arithm{\'e}tique (contr{\^o}le des d{\'e}nominateurs)
et des hypoth{\`e}ses de type analytique (uniformisation).

L'irruption de l'arithm{\'e}tique dans ces questions n'est pas nouvelle
mais remonte --- au moins --- {\`a} la communication~\cite{eisenstein1852}
d'Eisenstein o{\`u} ce dernier {\'e}tablit le th{\'e}or{\`e}me suivant:
\emph{soit $y=\sum_{n=0}^\infty a_n x^n\in\Q[[x]]$
le d{\'e}veloppement en s{\'e}rie de Taylor
en l'origine d'une fonction alg{\'e}brique. Alors, il existe un
entier $A \geq 1$ tel que, pour tout entier $n\geq 0$,
$a_n A^{n+1}$ est entier.} En particulier, les nombres premiers
qui divisent les d{\'e}nominateurs des rationnels $a_n$ sont en
nombre fini.
Une cons{\'e}quence imm{\'e}diate de ce r{\'e}sultat, cit{\'e}e d'ailleurs par Eisenstein,
est la transcendence des fonctions logarithme ou exponentielle,
{\og mais aussi de beaucoup d'autres.\fg}

\subsection{}

Comme il m'est impossible de donner un aper\c{c}u des m{\'e}thodes
et des r{\'e}sultats en th{\'e}orie des nombres transcendants,
je me vois oblig{\'e} de renvoyer le lecteur {\`a} l'immense
litt{\'e}rature sur le sujet, 
{\`a} laquelle le petit ouvrage collectif~\cite{bertrand-etal84}
peut {\^e}tre une bonne introduction, quoiqu'un peu ancienne.
Rappelons simplement que, classiquement,
une {\og d{\'e}monstration de transcendance\fg}
met successivement en jeu quatre ingr{\'e}dients:
on raisonne par l'absurde, puis 
\begin{enumerate}
\item on construit
une \emph{fonction auxiliaire}
{\`a} l'aide du \emph{lemme de Siegel}: celui-ci
fournit {\`a} un syst{\`e}me sous-d{\'e}termin{\'e}
d'{\'e}quations lin{\'e}aires {\`a} coefficients entiers une solution
enti{\`e}re de petite taille, mais \emph{non nulle}.
L'\emph{{\'e}valuation} de cette fonction auxiliaire fournit
un entier $\xi$ sur lequel va porter la contradiction ;

\item un certain nombre d'estimations de nature analytique,
telles que le \emph{lemme de Schwarz} et l'in{\'e}galit{\'e} d'Hadamard 
majorent $\abs{\xi}$ ; 

\item le fait (trivial) qu'un entier non nul est de valeur absolue
sup{\'e}rieure ou {\'e}gale {\`a}~$1$ implique, si la majoration
pr{\'e}c{\'e}dente est assez bonne, que $\xi=0$.
Dans les corps de nombres, ce principe est remplac{\'e}
par la \emph{formule du produit} ;

\item un \emph{lemme de z{\'e}ros}, de nature g{\'e}n{\'e}ralement alg{\'e}bro-g{\'e}om{\'e}trique,
montre qu'alors la fonction auxiliaire initiale est nulle,
d'o{\`u} la contradiction.
\end{enumerate}

Cependant, la m{\'e}thode des \emph{d{\'e}terminants d'interpolation}
introduite par M.~Laurent, cf.~\cite{laurent91},
ne fait pas intervenir le lemme de Siegel.
En fait, plut{\^o}t que sur une solution du syst{\`e}me lin{\'e}aire,
on peut raisonner directement sur les mineurs de sa matrice.

\subsection{}

La \emph{m{\'e}thode des pentes}, introduite par Bost {\`a} propos
d'un th{\'e}or{\`e}me de Masser et W{\"u}stholz (voir l'expos{\'e}~\cite{bost96}
{\`a} ce s{\'e}minaire) est une version g{\'e}om{\'e}trique {\og intrins{\`e}que\fg}
des d{\'e}terminants d'interpolation.
Sa formulation n{\'e}cessite le langage de la g{\'e}om{\'e}trie d'Arakelov.

Pla\c{c}ons-nous pour l'instant dans le cas le plus simple
o{\`u} les objets sont des r{\'e}seaux euclidiens $\bar E=(E,q_E)$, o{\`u}
$E$ est un $\Z$-module libre de rang fini et $q_E$ une forme quadratique 
d{\'e}finie positive sur $E_\R=E\otimes_\Z\R$.
{\`A} un tel objet, on peut associer un nombre r{\'e}el $\hdeg\bar E$,
son \emph{degr{\'e} arithm{\'e}tique}, qui n'est autre que l'oppos{\'e} du logarithme
du covolume de $E$ dans $E_\R$. En d'autres termes,
$(e_1,\dots,e_d)$ est une $\Z$-base de $E$,
\begin{equation}
\label{eq.gram}
 \hdeg \bar E =
     - \frac12 \log \det \big( q_E(e_i,e_j) \big)_{1\leq i,j\leq d}.
\end{equation}
On d{\'e}finit ensuite sa \emph{pente} qui est tout simplement
son degr{\'e} divis{\'e} par son rang: 
\[ \mu(\bar E) = \hdeg \bar E / \rg E. \]
Dans ce contexte, on a aussi une notion de polygone de Harder-Narasimhan
(cf.~\cite{stuhler76,grayson84,bost96})
mais ne nous int{\'e}ressera ici que sa \emph{plus grande pente},
$\mu_{\max}(\bar E)$:
c'est le maximum des pentes des sous-r{\'e}seaux 
de $\bar E$, c'est-{\`a}-dire des r{\'e}seaux $\bar F = (F,q_E|_F)$
o{\`u} $F$ est un sous-$\Z$-module non nul de~$E$. 

{\'E}tant donn{\'e}s deux tels r{\'e}seaux euclidiens $\bar E$ et $\bar F$,
un homomorphisme \mbox{$\phi\colon E\ra F$} poss{\`e}de une \emph{hauteur},
le logarithme de la norme de l'homomorphisme d'espaces euclidiens
induit:
\[ h(\phi) = \log \norm{\phi} = \frac12\log \sup_{e\in E_\R\setminus\{0\}}
        \frac{q_F(\phi(e))}{q_E(e)}. \]

Si $\phi$ est \emph{injectif}, on a alors une in{\'e}galit{\'e},
dite \emph{in{\'e}galit{\'e} de pentes}:
\begin{equation}\label{eq.inegalite}
 \mu (\bar E) \leq \mu_{\max}(\bar F)
           + h(\phi). 
\end{equation}
C'est essentiellement une reformulation de l'in{\'e}galit{\'e} d'Hadamard.
Le m{\'e}rite de cette in{\'e}galit{\'e} est de synth{\'e}tiser les diff{\'e}rentes {\'e}tapes
d'une preuve de transcendance: si le lemme de Siegel dispara{\^\i}t
(en apparence seulement, voir plus bas), les estim{\'e}es analytiques 
interviennent dans la majoration de $h(\phi)$. Quand l'in{\'e}galit{\'e}
obtenue est absurde, le morphisme $\phi$ ne peut pas {\^e}tre injectif,
et l{\`a} intervient {\'e}ventuellement le lemme de z{\'e}ros.

\subsection{}
Expliquons maintenant pourquoi, joint au th{\'e}or{\`e}me de Minkowski,
le lemme de Siegel est, {\`a} des facteurs num{\'e}riques pr{\`e}s,
un corollaire de l'in{\'e}galit{\'e} de pentes~\eqref{eq.inegalite}.
\emph{Soit donc $\Phi=(a_{ij})$ une matrice {\`a} coefficients entiers
{\`a} $r$ lignes et $n$ colonnes, avec $n>r$.
Si $A=\max \abs{a_{i,j}}$, le lemme de Siegel classique
affirme qu'il existe un {\'e}l{\'e}ment non nul
$\mathbf x= \transp (x_1,\dots,x_n)\in\Z^n$
tel que $\Phi\cdot\mathbf x=0$ et $\abs{x_i}\leq (nA)^{r/(n-r)}$
pour tout $i\in\{1;\dots;n\}$.}

Soit alors $\bar E$ le r{\'e}seau $\Z^n$ muni de la forme quadratique 
$q_E(x_1,\dots,x_n)=\sum_{i=1}^n x_i^2$,
soit $\bar F$ le r{\'e}seau analogue de rang~$r$,
et soit $\phi\colon E\ra F$ l'homomorphisme d{\'e}fini par
la matrice~$\Phi$.
Comme $n>r$, il n'est pas injectif et son noyau, muni
de la norme euclidienne induite, est un sous-r{\'e}seau
satur{\'e} $\bar E_1$. On d{\'e}finit alors $\bar E_2=\bar E/\bar E_1$
le r{\'e}seau quotient, muni de la norme euclidienne quotient.
Alors, $\phi$ induit un homomorphisme injectif 
$\phi_2\colon \bar E_2\ra \bar F$, si bien que l'in{\'e}galit{\'e}
de pentes~\eqref{eq.inegalite} s'{\'e}crit
\begin{equation}
 \mu(\bar E_2) \leq \mu_{\max}(\bar F) + h(\phi_2). 
\end{equation}
Par construction,
\begin{equation}
 h(\phi_2) = h(\phi) = \frac12\log \sup_{\mathbf x\in E_\R\setminus\{0\}}
   \frac{q_F(\phi(\mathbf x))}{q_E(\mathbf x)} 
        \leq \log (\sqrt r A).
\end{equation}
De plus, les covolumes de r{\'e}seaux sont multiplicatifs dans les suites
exactes, donc les degr{\'e}s arithm{\'e}tiques sont additifs et
\begin{equation}
 \mu(\bar E_2) = \rg(E_2)^{-1} \hdeg E_2 = \rg(E_2)^{-1} (\hdeg \bar E - \hdeg \bar
E_1) = -\frac{\rg (E_1)}{\rg(E_2)} \mu( \bar E_1)  
\end{equation}
car $\hdeg\bar E = r \hdeg (\Z,\norm{\cdot}) = 0$.
De m{\^e}me, $\hdeg\bar F=0$ et, 
plus g{\'e}n{\'e}ralement, le degr{\'e}
de tout sous-r{\'e}seau de $\bar F$ est n{\'e}gatif ou nul:
d'apr{\`e}s la formule de Gram (cf.~la formule~\eqref{eq.gram}),
le carr{\'e} du covolume d'un tel sous-r{\'e}seau de $\Z^n$ est en effet un entier
non nul, si bien que $\mu_{\max}(\bar F)=0$.
Ainsi, utilisant que $d=\rg E_1\geq n-r$, on a 
\begin{equation}
\label{eq.siegel}
 \mu(\bar E_1) \geq - \frac{r}{n-r} \log (A\sqrt r) 
\end{equation}

\'Ecrivons maintenant le th{\'e}or{\`e}me de Minkowski: il affirme
que \emph{si le volume (dans $E_{1,\R})$ de la boule $B_t$ de rayon $t$
est sup{\'e}rieur ou {\'e}gal {\`a} $2^d$ fois le covolume de $E_1$,
alors $B_t$ contient au moins un point non nul de $E_1$.}
Notons $\beta_d = \pi^{d/2}/\Gamma(d/2)$ le volume de la boule
unit{\'e} euclidienne en dimension $d$.
Retranscrit en termes de pentes, le th{\'e}or{\`e}me de Minkowski
affirme l'existence d'un point non nul $\mathbf x\in E_1$
tel que $q_E(\mathbf x)\leq t^2$ d{\`e}s que 
\begin{equation}
\label{eq.minkowski}
\mu(\bar E_1) \geq - \log t + \log (2\beta_d^{-1/d}). 
\end{equation}

Ainsi, il existe une solution non nulle $\mathbf x=\transp (x_1,\dots,x_n)\in\Z^n$
au syst{\`e}me $\Phi\cdot\mathbf x=0$ v{\'e}rifiant
\begin{equation}
 \norm{\mathbf x}= \big(\sum_{i=1}^n x_i^2\big)^{1/2} \leq 
         \big(A\sqrt r\big)^{r/(n-r)}  2 \beta_d^{-1/d}.
\end{equation}
De plus, la formule de Stirling montre que lorsque $d$ tend
vers~$+\infty$,
\[ 2\beta_d^{-1/d} \simeq \sqrt{2d/e}. \]

\subsection{}
Dans le m{\^e}me volume que celui o{\`u} est publi{\'e}~\cite{chudnovsky85b},
les Chudnovsky utilisent des techniques similaires
pour d{\'e}buter la th{\'e}orie arithm{\'e}tique des $G$-fonctions de Siegel,
introduites dans~\cite{siegel29}.
Faute de place, nous ne dirons rien de cette th{\'e}orie
qui pourtant a vu r{\'e}cemment quelques d{\'e}veloppement majeurs,
suite notamment aux travaux de Bombieri~\cite{bombieri81}
et~Andr{\'e}. Citons ainsi une th{\'e}orie
g{\'e}om{\'e}trique des op{\'e}rateurs diff{\'e}rentiels de type $G$
(Andr{\'e}, Baldassarri~\cite{andre-b97}) ainsi que la g{\'e}n{\'e}ralisation
en dimension sup{\'e}rieure,
due {\`a} L.~Di Vizio~\cite{divizio2000},
du th{\'e}or{\`e}me des Chudnovsky affirmant qu'un op{\'e}rateur diff{\'e}rentiel
minimal annulant une $G$-fonction est un $G$-op{\'e}rateur.
Enfin, signalons deux articles r{\'e}cents
d'Andr{\'e}~\cite{andre2000a,andre2000b}
consacr{\'e}s {\`a} une {\og th{\'e}orie Gevrey arithm{\'e}tique\fg}
dont le th{\'e}or{\`e}me de Chudnovsky {\'e}voqu{\'e} est un outil essentiel.
Il en d{\'e}duit de remarquables applications, par exemple
une nouvelle preuve du th{\'e}or{\`e}me de Siegel-Shidlovsky
et, inspir{\'e} par la preuve de B{\'e}zivin-Robba~\cite{bezivin-r89},
deux (\string!) du th{\'e}or{\`e}me de Lindemann-Weierstra{\ss}.

\bigskip
Je remercie Y.~Andr{\'e}, D.~Bertrand, J.-B.~Bost et Y.~Laszlo
pour l'aide qu'ils m'ont apport{\'e}e pendant la pr{\'e}paration
de cet expos{\'e}. Je remercie aussi C.~Gasbarri et A.~Thuillier 
pour leur lecture attentive.

\section{Quatre th{\'e}or{\`e}mes}\label{sec.4thms}

\subsection{}
Les {\'e}nonc{\'e}s que nous pr{\'e}sentons maintenant font intervenir
des objets alg{\'e}bro-g{\'e}om{\'e}triques d{\'e}finis sur un corps de nombres
(nombres, vari{\'e}t{\'e}s alg{\'e}briques,
{\'e}quations diff{\'e}rentielles,  sous-alg{\`e}bres de Lie, etc.)
v{\'e}rifiant une propri{\'e}t{\'e} modulo $p$ pour presque tout nombre
premier $p$.
Le sens de ce genre d'hypoth{\`e}ses est le suivant.
Soit $K$ un corps de nombres et soit $\mathfrak o_K$
son anneau d'entiers.
Un objet alg{\'e}bro-g{\'e}om{\'e}trique sur $K$ {\og de type fini\fg}
est d{\'e}fini par une famille finie d'{\'e}l{\'e}ments de $K$ (par exemple,
dans le cas d'une vari{\'e}t{\'e} alg{\'e}brique,
les coefficients des {\'e}quations polynomiales qui la d{\'e}finissent) ; ces
{\'e}l{\'e}ments ont un d{\'e}nominateur commun $D\in\Z$, si bien qu'en fait
l'objet initial {\og vit\fg} naturellement sur l'anneau
$\mathfrak o_K[1/D]$. Il y a bien s{\^u}r des choix mais
{\'e}tant donn{\'e}es deux structures enti{\`e}res sur $\mathfrak o_K[1/D]$
et $\mathfrak o_K[1/D']$, il existera un multiple commun {\`a} $D$ et $D'$,
soit $D''$, tel que les structures enti{\`e}res co{\"\i}ncident, une fois
regard{\'e}es sur $\mathfrak o_K[1/D'']$.
C'est bien s{\^u}r dans le langage des \emph{sch{\'e}mas} que de telles
consid{\'e}rations trouvent leur place naturelle.

Les id{\'e}aux maximaux d'un tel anneau $\mathfrak o_K[1/D]$
s'identifient naturellement {\`a} une partie de compl{\'e}mentaire fini
de l'ensemble des id{\'e}aux maximaux de $\mathfrak o_K$ (ceux qui
ne contiennent pas $D$). R{\'e}duire modulo $p$ pour presque tout $p$
signifie que, pour tout id{\'e}al maximal $\mathfrak p$ de $\mathfrak o_K[1/D]$
(sauf peut-{\^e}tre un nombre fini d'entre eux),
on consid{\`e}re la structure enti{\`e}re modulo l'id{\'e}al $\mathfrak p$.
On obtient ainsi une structure alg{\'e}brique analogue {\`a}
celle de d{\'e}part, mais sur un corps fini.
C'est sur ces structures {\og modulo $p$\fg} que portent les hypoth{\`e}ses
des th{\'e}or{\`e}mes.

Nous donnerons volontairement les {\'e}nonc{\'e}s dans ce style informel,
en en indiquant juste apr{\`e}s la signification pr{\'e}cise.

\begin{Theo}\label{theo.kronecker}
Soit $\alpha$ un nombre alg{\'e}brique tel que pour presque
tout nombre premier $p$, la r{\'e}duction modulo~$p$ de $\alpha$ 
est dans le sous-corps premier $\F_p$.
Alors, $\alpha$ est un nombre rationnel.
\end{Theo}
Autrement dit, $\alpha$ est un {\'e}l{\'e}ment d'un corps de nombres $K$.
{\'E}crivons le $\beta/D$ o{\`u} $D$ est un entier $\geq 1$
et $\beta$ un entier alg{\'e}brique.
L'hypoth{\`e}se signifie que pour presque tout id{\'e}al maximal $\mathfrak p$
de $\mathfrak o_K$, il existe un {\'e}l{\'e}ment $n_{\mathfrak p}\in\Z$
tel que $\beta\equiv n_{\mathfrak p}D\pmod{\mathfrak p}$.
La conclusion est alors que $\beta\in\Z$.

Ce r{\'e}sultat appara{\^\i}t dans l'article~\cite{kronecker1880} de Kronecker
qui le d{\'e}rive du comportement
du logarithme de la fonction z{\^e}ta de Dedekind du corps $\Q(\alpha)$ en
$s=1$.\footnote{L'argument de Kronecker semble cependant incomplet.
L'existence d'une densit{\'e} de certains ensembles de nombres premiers,
aujourd'hui cons{\'e}quence du th{\'e}or{\`e}me de \v Cebotarev, est en
effet utilis{\'e}e sans justification.}
Aujourd'hui, il appara{\^\i}t souvent comme une cons{\'e}quence 
du th{\'e}or{\`e}me de densit{\'e} de \v Cebotarev.
En lien avec le th{\'e}or{\`e}me~\ref{theo.diff-exact} ci-dessous,
les Chudnovsky en fournissent une d{\'e}monstration 
{\og diophantienne\fg} dans~\cite{chudnovsky85b}, 
c'est-{\`a}-dire fond{\'e}e sur les m{\'e}thodes introduites par la
th{\'e}orie des nombres transcendants.

\begin{Theo}\label{theo.tate-ell}
Soit $E$ et $E'$ deux courbes elliptiques sur $\Q$.
Alors, $E$ et $E'$ sont isog{\`e}nes sur $\Q$
si et seulement si, pour presque tout nombre premier $p$,
$E$ et $E'$ ont le m{\^e}me nombre de points sur $\F_p$.
\end{Theo}
Donnons l{\`a} encore une interpr{\'e}tation na{\"\i}ve:
partant de deux courbes elliptiques $E$ et $E'$ sur $\Q$
on peut en trouver des {\'e}quations
(affines, sous forme de Weierstra\ss) $y^2=4x^3+ax+b$ et $y^2=4x^3+a'x+b'$,
avec $a$, $b$, $a'$, $b'$ dans $\Z$.
Pour tout nombre premier $p$ sauf $2$, $3$
et ceux qui divisent le produit
des discriminants $4a^3+27b^2$ et $4a^{\prime 3}+27b^{\prime 2}$,
ces {\'e}quations d{\'e}finissent des courbes elliptiques $E_p$
et $E'_p$ sur le corps fini $\F_p$.

Une $\Q$-isog{\'e}nie entre $E$ et $E'$, c'est-{\`a}-dire en l'occurrence
un morphisme non constant d{\'e}fini sur $\Q$ fournira pour tous ces nombres
premiers (sauf quelques-uns, peut-{\^e}tre) une isog{\'e}nie entre
les courbes elliptiques $E_p$ et $E'_p$. Il est alors connu
que $E_p$ et $E'_p$ ont m{\^e}me nombre de points sur $\F_p$.
La r{\'e}ciproque est le point difficile.

D{\'e}montr{\'e} par Serre~\cite{serre68b}
si l'un des invariants $j(E)$ et $j(E')$ n'est pas entier,
c'est un cas particulier du {\og th{\'e}or{\`e}me d'isog{\'e}nie
de Faltings\fg}, cf.~\cite{faltings83}, 
lequel fournit un {\'e}nonc{\'e} analogue pour deux vari{\'e}t{\'e}s
ab{\'e}liennes d{\'e}finies sur un corps de nombres.
Dans~\cite{chudnovsky85b},
les Chudnovsky avaient donn{\'e} une d{\'e}monstration
diophantienne du th{\'e}or{\`e}me~\ref{theo.tate-ell}.
Cependant, leur d{\'e}monstration n{\'e}cessitait en outre
le th{\'e}or{\`e}me de Cartier et Honda~\cite{honda70}
reliant la fonction z{\^e}ta de Hasse-Weil d'une $\Q$-courbe
elliptique et le groupe formel de son mod{\`e}le de N{\'e}ron sur $\Z$.
Dans~\cite{graftieaux2001},
Graftieaux a donn{\'e} une d{\'e}monstration analogue du th{\'e}or{\`e}me d'isog{\'e}nie pour
les vari{\'e}t{\'e}s ab{\'e}liennes {\`a} multiplication r{\'e}elle d{\'e}finies sur $\Q$.
Cette preuve n{\'e}cessite une g{\'e}n{\'e}ralisation du th{\'e}or{\`e}me de Cartier 
et Honda, d{\'e}montr{\'e}e par Deninger et Nart dans~\cite{deninger-n90},
qui fournit un isomorphisme
entre les groupes formels des vari{\'e}t{\'e}s ab{\'e}liennes
en question. Graftieaux d{\'e}montre ensuite l'alg{\'e}bricit{\'e} de cet isomorphisme,
autrement dit l'alg{\'e}bricit{\'e} du graphe, lequel est un sous-groupe
formel du produit des deux vari{\'e}t{\'e}s ab{\'e}liennes.
Graftieaux a ensuite {\'e}tabli dans~\cite{graftieaux2001b}
un crit{\`e}re g{\'e}n{\'e}ral d'alg{\'e}bricit{\'e} de sous-groupes
formels d'une vari{\'e}t{\'e} ab{\'e}lienne d{\'e}finie sur un corps de nombres.
Ses r{\'e}sultats sont de plus effectifs: les hypoth{\`e}ses ne font intervenir
qu'un nombre fini explicite de nombres premiers.

Plus g{\'e}n{\'e}ralement, on a le th{\'e}or{\`e}me suivant, d{\^u} {\`a} Bost.

\begin{Theo}\label{theo.lie}
Soit $G$ un groupe alg{\'e}brique d{\'e}fini sur un corps
de nombres~$K$, soit $\mathfrak g$ son alg{\`e}bre de Lie.
Soit $\mathfrak h$ une sous-$K$-alg{\`e}bre de Lie de $\mathfrak g$
v{\'e}rifiant la propri{\'e}t{\'e} suivante:
pour presque tout nombre premier $p$, la r{\'e}duction
de $\mathfrak h$ modulo $p$
est une $p$-alg{\`e}bre de Lie.
Alors, il existe un sous-groupe alg{\'e}brique $H$ de $G$,
d{\'e}fini sur $K$, dont $\mathfrak h$ est l'alg{\`e}bre de Lie.
\end{Theo}

L'ingr{\'e}dient nouveau est la notion de $p$-alg{\`e}bre de Lie.
Si $G$ est un groupe alg{\'e}brique lisse sur un corps $k$,
son alg{\`e}bre de Lie $\mathfrak g$ est par d{\'e}finition le $k$-espace
vectoriel des d{\'e}rivations invariantes par translations sur $G$,
muni du crochet
$[D_1,D_2]=D_1D_2-D_2D_1$. En effet, si 
$D_1$ et $D_2$
sont deux d{\'e}rivations invariantes, leur commutateur en est encore une.
De plus, si $k$ est de caract{\'e}ristique $p>0$
et si $D$ est une d{\'e}rivation invariante, la formule du bin{\^o}me
montre que, pour tous germes de fonctions $f$ et $g$, on a
\[ D^p(fg) = \sum_{i=0}^p \binom pi D^i(f) D^{p-i}(g)
     = g D^p(f) + fD^p(g), \]
si bien que $D^p$ est encore une d{\'e}rivation.
Cette op{\'e}ration de puissance~\pieme
donne naissance {\`a} la notion de $p$-alg{\`e}bre de Lie,
cf.~\cite{jacobson79}, 5.7, \cite{borel91}, 1.3
ainsi que l'expos{\'e}~\cite{mathieu2000} {\`a} ce s{\'e}minaire.
{\'E}tant donn{\'e}e une telle $p$-alg{\`e}bre de Lie $\mathfrak g$, une
sous-$p$-alg{\`e}bre de Lie est alors une sous-alg{\`e}bre de Lie $\mathfrak h$
telle que pour tout $D\in\mathfrak h$, $D^p\in\mathfrak h$.

L'explicitation du {\og modulo~$p$\fg} se fait alors comme pr{\'e}c{\'e}demment.
Si $\mathfrak o_K$ est l'anneau des entiers de $K$,
il existe un entier $N\geq 1$ et un sch{\'e}ma en groupes
lisse $\mathscr G$ sur $\Spec\mathfrak o_K[1/N]$ dont
la fibre g{\'e}n{\'e}rique $\mathscr G\otimes K $ est $G$.
Son alg{\`e}bre de Lie $\Lie(\mathscr G)$
est une $\mathfrak o_K[1/N]$-alg{\`e}bre de Lie
telle que $\Lie(\mathscr G)\otimes K=\mathfrak g$.
Pour tout id{\'e}al maximal $\mathfrak p $ de $\mathfrak o_K$
ne contenant pas $N$, on dispose alors
d'un groupe alg{\'e}brique lisse sur le corps fini
$\F_{\mathfrak p}=\mathfrak o_K/\mathfrak p$ dont l'alg{\`e}bre de Lie
n'est autre que $\Lie(\mathscr G)\otimes \F_{\mathfrak p}$.
Si l'on note $p$ la caract{\'e}ristique du corps~$F_{\mathfrak p}$,
c'est en particulier une $p$-alg{\`e}bre de Lie.

Quitte {\`a} remplacer l'entier $N$ par un multiple,
une sous-alg{\`e}bre de Lie $\mathfrak h$ de $\mathfrak g$
d{\'e}finit de m{\^e}me une sous-alg{\`e}bre de Lie de $\Lie(\mathscr G)$,
et par cons{\'e}quent, par r{\'e}duction modulo~$\mathfrak p$,
des sous-alg{\`e}bres de Lie de $\Lie(\mathscr G)\otimes \F_{\mathfrak p}$.
L'hypoth{\`e}se est donc que pour presque tout id{\'e}al
maximal $\mathfrak p$, ces sous-alg{\`e}bres de Lie sont 
des sous-$p$-alg{\`e}bres de Lie.

\subsection{}
Le th{\'e}or{\`e}me~\ref{theo.lie} est d{\'e}montr{\'e} par Bost dans~\cite{bost2001}.
Il avait {\'e}t{\'e} ind{\'e}pendamment conjectur{\'e} par Ekedahl et Shepherd-Barron
dans leur article~\cite{ekedahl-sb2000}.
Il n'est peut-{\^e}tre pas inutile d'expliquer en quoi
les deux th{\'e}or{\`e}mes~\ref{theo.kronecker} et~\ref{theo.tate-ell}
en sont des cas particuliers.

\subsubsection{Th{\'e}or{\`e}me de Kronecker}
Consid{\'e}rons le groupe alg{\'e}brique $G_k=\Gm\times\Gm$
sur un corps $k$.
On note $x$ et $y$ les coordonn{\'e}es sur les deux produits, l'{\'e}l{\'e}ment
neutre {\'e}tant $(1,1)$. Alors, une base de l'alg{\`e}bre
de Lie $\mathfrak g_k$ de $G_k$
est constitu{\'e}e des champs de vecteurs $x\frac\partial{\partial x}$
et $y\frac{\partial}{\partial y}$.
et $\mathfrak g_k$ est la $k$-alg{\`e}bre de Lie commutative
$k^2$.

Si $k$ est de caract{\'e}ristique positive $p$, d{\'e}terminons
l'op{\'e}ration
de puissance~\pieme sur $\mathfrak g_k$. Il suffit
de la calculer dans le cas de $\Gm$, auquel cas on a
\begin{equation}
 \big(x\frac{\partial}{\partial x}\big)^p
         = x\frac{\partial}{\partial x}. 
\end{equation}
En effet, $x\frac{\partial}{\partial x}$ est l'unique d{\'e}rivation
invariante sur $\Gm$ qui envoie sur elle-m{\^e}me la fonction r{\'e}guli{\`e}re~$x$.
Ainsi, utilisant le fait que l'alg{\`e}bre de Lie de $(\Gm)^2$ est
commutative,
\[ \big(a x\frac\partial{\partial x}+b y \frac\partial{\partial y}\big)^p
     =a^p x \frac\partial{\partial x}+b^p y \frac\partial{\partial y}
.\]

Pour {\'e}tablir le th{\'e}or{\`e}me~\ref{theo.kronecker},
consid{\'e}rons un {\'e}l{\'e}ment $\alpha$ d'un corps de nombres $K$
et soit $G$ le groupe $\Gm\times\Gm$ sur $K$, d'alg{\`e}bre
de Lie $\mathfrak g=K^2$.
La droite $\mathfrak h=(1,\alpha)K$ est une sous-alg{\`e}bre de Lie.
Soit $\mathfrak p$ un id{\'e}al maximal de $\mathfrak o_K$
tel que $\alpha$ soit $\mathfrak p$-entier.
La r{\'e}duction modulo $\mathfrak p$ de $\mathfrak h$
est la sous-alg{\`e}bre de Lie de $\mathfrak h_{\mathfrak p}$
engendr{\'e}e par $(1,\alpha\mod\mathfrak p)$ dans $(\mathfrak o_K/\mathfrak p)^2$.
C'est une sous-$p$-alg{\`e}bre de Lie 
si et seulement si $(1^p,(\alpha \mod\mathfrak p)^p)$ est colin{\'e}aire {\`a} 
$(1,\alpha\mod\mathfrak p)$, c'est-{\`a}-dire si et seulement
si $\alpha^p= \alpha\pmod{\mathfrak p}$, ou encore si
$\alpha\mod\mathfrak p$ est dans le corps premier $\F_p$.
Par suite, sous l'hypoth{\`e}se du th{\'e}or{\`e}me de Kronecker,
l'alg{\`e}bre de Lie $\mathfrak h$ est l'alg{\`e}bre de Lie d'un
$K$-sous-groupe alg{\'e}brique $H$ de $\Gm\times\Gm$. Or, ceux-ci sont
d{\'e}finis par une {\'e}quation de la forme $x^m=y^n$ pour deux entiers $m$ et $n$.
L'alg{\`e}bre de Lie d'un tel $H$ est la droite d'{\'e}quation
$m\xi=n\eta$ dans le plan $K^2$ de coordonn{\'e}es $(\xi,\eta)$,
donc est engendr{\'e}e par le vecteur $(n,m)$.
On a ainsi $\alpha=m/n$. 

\subsubsection{Th{\'e}or{\`e}me d'isog{\'e}nie}
Pour commencer, soit $E$ une courbe elliptique sur un
corps $k$.
Son alg{\`e}bre de Lie est une $k$-alg{\`e}bre de Lie de dimension 1,
bien entendu commutative. De plus, on a une identification
canonique $\Lie(E)=H^1(E,\mathscr O_E)$.
Si $k$ est de caract{\'e}ristique positive $p$, l'op{\'e}ration de puissance
\pieme sur $\Lie(E)$ se confond alors
avec l'action du morphisme dit {\og Frobenius absolu\fg}
$F\colon E\ra E$. 
Supposons de plus que $k$ est le corps premier $\F_p$.
Alors, $F$ induit un endomorphisme $p$-lin{\'e}aire, donc lin{\'e}aire,
de $H^1(E,\mathscr O_E)$ ; c'est ainsi une homoth{\'e}tie 
dont nous noterons $A\in\F_p$ le rapport (\emph{invariant de Hasse}).
On sait d'autre part depuis Hasse qu'il existe deux entiers alg{\'e}briques
$\alpha$ et $\beta$ tels que $\# E(\F_p)=p+1-(\alpha+\beta)$ et
$\alpha\beta=p$. De plus, modulo~$p$, l'un des deux, disons $\alpha$,
s'interpr{\`e}te comme l'action de l'endomorphisme~$F$
sur $\Lie E$, tandis que $\beta$ correspond {\`a} l'action du Verschiebung
(d{\'e}fini par $FV=p$),
c'est-{\`a}-dire~$A$. Par suite $\# E(\F_p)\equiv 1-A$ modulo~$p$.

Pla\c{c}ons-nous maintenant sous les hypoth{\`e}ses du th{\'e}or{\`e}me~\ref{theo.tate-ell}
et consid{\'e}rons le groupe alg{\'e}brique $G=E\times E'$ sur $\Q$,
d'alg{\`e}bre de Lie $\mathfrak g=\Lie(E)\oplus \Lie(E')$.
Soit $\mathfrak h\subset\mathfrak g$
une droite arbitraire se projetant surjectivement
sur les deux facteurs. L'hypoth{\`e}se que $E$ et $E'$ ont, pour
presque tout nombre premier $p$, m{\^e}me nombre de points
modulo~$p$ implique que,  modulo~$p$ pour presque tout $p$,
$\mathfrak h$ d{\'e}finit une sous-$p$-alg{\`e}bre de Lie de la r{\'e}duction
modulo~$p$ de $\mathfrak g$.
C'est ainsi l'alg{\`e}bre de Lie d'un sous-groupe alg{\'e}brique $H\subset E\times E'$
d{\'e}fini sur $\Q$.
Or, un tel sous-groupe fournit automatiquement une isog{\'e}nie
entre $E$ et $E'$ (consid{\'e}rer par exemple $H$ comme une correspondance).

\begin{Theo}\label{theo.diff-exact}
Soit $X$ une vari{\'e}t{\'e} alg{\'e}brique lisse d{\'e}finie
sur un corps de nombres et soit $\omega$ une forme diff{\'e}rentielle sur $X$.
On suppose que pour presque tout nombre premier $p$,
$\omega$ est modulo $p$ une forme diff{\'e}rentielle exacte ($\omega=\mathrm df$).
Alors,  $\omega$ est une forme diff{\'e}rentielle exacte.

On a un r{\'e}sultat similaire dans le cas logarithmique: si pour
presque tout nombre premier $p$, $\omega$ est modulo $p$
une forme logarithmique exacte ($\omega = \mathrm d\log f=\mathrm df/f$),
alors, il existe un entier $n\geq 1$ tel que $n\omega$
est une forme diff{\'e}rentielle logarithmiquement exacte.
\end{Theo}
Ce th{\'e}or{\`e}me est d{\'e}montr{\'e}
par Andr{\'e}~\cite{andre89} dans le premier cas
et
par les Chudnovsky
dans~\cite{chudnovsky85b} dans le cas logarithmique.
Il {\'e}tait conjectur{\'e} dans le premier cas par Ogus dans~\cite[\S 2]{ogus82}
et, en liaison avec la conjecture de Grothendieck,
par Katz dans~\cite[p.~2]{katz72} pour le cas logarithmique.

Comme exemple d'application, montrons
comment en d{\'e}duire le th{\'e}or{\`e}me~\ref{theo.kronecker}.
Consid{\'e}rons la courbe $X=\Gm=\Spec K[x,x^{-1}]$ sur un corps de nombres $K$,
et, si $\alpha\in K$, posons $\omega= \alpha x^{-1} \mathrm dx$.
Si $\mathfrak p$ est un id{\'e}al maximal de $\mathfrak o_K$
et $n_{\mathfrak p}\in\Z$ un entier tel que $\alpha - n_{\mathfrak
p}\in\mathfrak p$, $\omega $ modulo $\mathfrak p$ est {\'e}gale
{\`a} $n_{\mathfrak p} x^{-1}\mathrm dx$, donc est la diff{\'e}rentielle
logarithmique de $x^{n_ {\mathfrak p}}$. D'apr{\`e}s le
th{\'e}or{\`e}me~\ref{theo.diff-exact}, il existe $n\geq 1$ tel que
$n\alpha x^{-1}\mathrm dx$ est une forme diff{\'e}rentielle logarithmiquement
exacte, c'est-{\`a}-dire $n\alpha\in\Z$, d'o{\`u} $\alpha\in\Q$.

Sur la droite projective, le th{\'e}or{\`e}me~\ref{theo.diff-exact}
admet le corollaire suivant:
\emph{soit $y\in\Z[[x]]$ une s{\'e}rie formelle {\`a} coefficients
entiers dont la d{\'e}riv{\'e}e est une fonction alg{\'e}brique, alors
$y$ est une fonction alg{\'e}brique.}
La variante de cet {\'e}nonc{\'e} o{\`u} alg{\'e}brique est remplac{\'e}
par \emph{rationnelle} a {\'e}t{\'e} d{\'e}montr{\'e}e par P\'olya ;
c'est une application classique du th{\'e}or{\`e}me de Borel~\cite{eborel1894}.
Signalons d'ailleurs que ce crit{\`e}re a {\'e}t{\'e} g{\'e}n{\'e}ralis{\'e}
par B{\'e}zivin et Robba dans~\cite{bezivin-r89b} au cas d'op{\'e}rateurs
diff{\'e}rentiels d'ordre sup{\'e}rieur. Cette g{\'e}n{\'e}ralisation
leur a permis d'en d{\'e}duire une nouvelle d{\'e}monstration
du th{\'e}or{\`e}me de Lindemann-Weierstrass.

F.~Caligari en a donn{\'e} une application aux courbes modulaires:
joint au th{\'e}or{\`e}me de Manin-Drinfel'd,
il implique en effet qu'une forme modulaire
(non n{\'e}cessairement cuspidale)
de poids~$2$ qui pour presque tout nombre premier~$p$
est fix{\'e}e par l'op{\'e}rateur~$T_p$ modulo~$p$ est une combinaison
lin{\'e}aire de s{\'e}ries d'Eisenstein.

\section{La conjecture de Grothendieck}\label{sec.grothendieck}

\subsection{}
La \emph{conjecture de Grothendieck} est un crit{\`e}re arithm{\'e}tique
qui pr{\'e}dit qu'un syst{\`e}me diff{\'e}rentiel lin{\'e}aire, disons
de la forme
\begin{equation}\label{eq.sys-diff}
  \frac{\mathrm d}{\mathrm dz} Y = A(z) Y, \qquad A(z)\in M_d(\Q(z)) 
\end{equation}
poss{\`e}de une base de solutions alg{\'e}briques, c'est-{\`a}-dire
dont les solutions holomorphes au voisinage 
d'un point $z_0\in\Q$ qui n'est pas un p{\^o}le de $A$ sont alg{\'e}briques
sur $\Q(z)$.
La condition, conjecturalement suffisante, est
que pour presque tout nombre premier $p$, le syst{\`e}me
diff{\'e}rentiel 
obtenu par r{\'e}duction modulo $p$ ait une base de solutions
alg{\'e}briques sur $\F_p(z)$.

Le point fondamental, d{\^u} {\`a} Cartier et Honda, est que cette derni{\`e}re
condition est, pour $p$ un nombre premier fix{\'e}, effectivement v{\'e}rifiable. 
Pour rester {\'e}l{\'e}mentaire, d{\'e}finissons une suite $(A_n)$ de matrices
$n\times n$ {\`a} coefficients dans $\Q(z)$ par
\begin{equation}
 A_0(z) = I_d, \quad  A_1(z)=A(z), \quad
A_{n+1}(z) = \frac{\mathrm d}{\mathrm dz} A_n(z) +  A_n(z)A(z).
\end{equation}
Si l'on peut r{\'e}duire $A$ modulo $p$, alors on a l'{\'e}quivalence
des propositions suivantes:
\begin{itemize}
\item le syst{\`e}me diff{\'e}rentiel~\eqref{eq.sys-diff} modulo $p$ admet une base
de solutions alg{\'e}briques sur $\F_p(z)$ ;
\item il admet une base de solutions dans $\F_p(z)$ ;
\item il admet une base de solutions dans $\F_p((z))$ ;
\item la matrice de {\og $p$-courbure\fg} $A_p$ est nulle modulo $p$.
\end{itemize}
Sous cette forme, cet {\'e}nonc{\'e} est assez {\'e}l{\'e}mentaire. Si $Y(z)$ est une solution
du syst{\`e}me~\eqref{eq.sys-diff}, on v{\'e}rifie par r{\'e}currence que pour tout
$n$, $\frac{\mathrm d^n}{\mathrm dz^n} Y(z)=A_n(z) Y(z)$.
Si $Y$ est une solution de ce syst{\`e}me, par exemple 
dans $\F_p((z))$, on a
$(\mathrm d^p/\mathrm dz^p) Y(z)=0$. L'existence d'une
\emph{base} de solutions dans $\F_p((z))$ implique alors que $A_p(z)=0$
modulo $p$.
Dans l'autre sens, supposant pour simplifier que~$0$
n'est pas une singularit{\'e} de~$A$,  on constate que la formule 
(de Taylor !) 
\begin{equation}
 Y(z) =  \big (\sum_{i=0}^{p-1} \frac{(-z)^i}{i!} A_i(z) \big) ^{-1}  
\end{equation}
fournit une matrice fondamentale de solutions dans $\F_p(z)$.

\subsection{}
Ces consid{\'e}rations s'{\'e}tendent {\`a} un contexte plus g{\'e}n{\'e}ral que
nous pr{\'e}sentons maintenant.
Soit $X$ une vari{\'e}t{\'e} alg{\'e}brique lisse sur un corps de nombres $K$
et soit $E$ un $\mathscr O_X$-module localement libre de rang fini.
Une connexion sur $E$ est une application\break $K$-lin{\'e}aire
\begin{equation}
\nabla\colon E\ra E\otimes_{\mathscr O_X} \Omega^1_X
\end{equation}
v{\'e}rifiant $\nabla (f e) = f \nabla(e) + e\otimes \mathrm df$
pour toutes sections locales $f$ de $\mathscr O_X$
et $e$ de $E$.
Il sera pratique de consid{\'e}rer l'homomorphisme dual
\begin{equation}\label{eq.tx-end}
 \mathrm TX \ra \CEnd(E), \qquad  \partial \mapsto \nabla(\partial),
\end{equation}
$\mathrm TX$ d{\'e}signant le fibr{\'e} tangent de $X$.
On notera aussi $E^\nabla$ le faisceau des \emph{sections
horizontales} de $E$, c'est-{\`a}-dire des germes
de sections locales $e$ de $E$ telles  que $\nabla(e)=0$.
Par d{\'e}finition, $E$ est \emph{trivial} s'il est engendr{\'e}
par $E^\nabla$ comme $\mathscr O_X$-module.
On dit aussi que la connexion $\nabla$ est \emph{int{\'e}grable}
si l'homomorphisme~\eqref{eq.tx-end} est un homomorphisme
d'alg{\`e}bres de Lie, c'est-{\`a}-dire si
pour tous champs de vecteurs locaux $\partial_1$ et $\partial_2$,
\begin{equation}
[\nabla(\partial_1),\nabla(\partial_2)] = \nabla([\partial_1,\partial_2]).
\end{equation}
% Une connexion $\nabla$ sur $E$ s'{\'e}tend par la formule de Leibniz
% en une famille de connexions 
% \[\nabla \colon E\otimes \Omega^i_X \ra E\otimes\Omega^{i+1}_X,
%  \qquad e\otimes\omega \mapsto \nabla(e)\wedge \omega + e\otimes \mathrm
% d\omega. \]
% De plus, un calcul montre que $\nabla^2$ est
% $\mathscr O_X$-lin{\'e}aire et qu'il existe une
% $2$-forme $\psi(E,\nabla)\in\Gamma(X,\CEnd(E)\otimes_{\mathscr
% O_X}\Omega^2_X)$, la \emph{courbure} de $(E,\nabla)$,
% telle que $\nabla^2=\psi(E,\nabla)\wedge\cdot$.
% si sa courbure est nulle.
% Cela signifie aussi que $\nabla$ commute au crochet de Lie:
C'est bien s{\^u}r {\'e}quivalent au fait que la \emph{courbure}
$\psi(E,\nabla)=\nabla^2\in\End(E)\otimes\Omega^2_X$ soit nulle.

Supposons $\nabla$ int{\'e}grable.
Si $K$ est un corps de caract{\'e}ristique positive $p$,
$\mathrm TX$ et $\CEnd(E)$ sont des $p$-alg{\`e}bres de Lie et on
d{\'e}finit une \emph{$p$-courbure} qui mesure le d{\'e}faut
pour l'homomorphisme~\eqref{eq.tx-end} d'{\^e}tre un homomorphisme
de $p$-alg{\`e}bres de Lie:
\begin{equation}
\psi_p \colon \mathrm T_X \ra \CEnd(E), \qquad
       \partial \mapsto \nabla(\partial)^p -\nabla(\partial^p).
\end{equation}
C'est une application additive, $p$-lin{\'e}aire:
\[ \psi_p( f_1 \partial_1 + f_2\partial_2)
        = f_1^p \psi_p(\partial_1) + f_2^p \psi_p(\partial_2). \]
Le th{\'e}or{\`e}me de Cartier affirme alors l'{\'e}quivalence des propri{\'e}t{\'e}s
suivantes (cf.~\cite{katz70}, th{\'e}or{\`e}me~5.1):
\begin{itemize}
\item pour toute section locale $\partial $ de $\mathrm TX$,
$\psi_p(\partial)=0$ ;
\item  l'homomorphisme~\eqref{eq.tx-end} est un homomorphisme
de $p$-alg{\`e}bres de Lie ;
\item le faisceau $E$ est trivial.
\end{itemize}

\begin{Conj}[Grothendieck]
Soit $X$ une vari{\'e}t{\'e} lisse sur un corps de nombres $K$
et $(E,\nabla)$ un module {\`a} connexion int{\'e}grable sur $X$.
On suppose que pour presque tout nombre premier $p$,
la r{\'e}duction modulo $p$ de $(E,\nabla)$ est {\`a} $p$-courbure nulle.
Alors, il existe un rev{\^e}tement {\'e}tale $f\colon Y\ra X$
tel que $f^*E$ soit trivial.
\end{Conj}

La condition pour un module {\`a} connexion int{\'e}grable
$(E,\nabla)$ d'avoir (presque toutes) ses $p$-courbures nulles a {\'e}t{\'e}
{\'e}tudi{\'e}e de mani{\`e}re approfondie par Katz.
En particulier, il est {\'e}tabli dans~\cite{katz70} qu'alors,
$(E,\nabla)$ est {\`a} singularit{\'e}s r{\'e}guli{\`e}res et ses exposants
sont des nombres rationnels, ce qui n'est d'ailleurs pas
sans rapport avec le th{\'e}or{\`e}me~\ref{theo.kronecker}
(voir aussi~\cite{dwork-g-s94}, III.6.1).

Il y a eu essentiellement deux approches de cette conjecture,
l'une g{\'e}om{\'e}trique, l'autre arithm{\'e}tique.
Rappelons que les p{\'e}riodes d'une famille de vari{\'e}t{\'e}s alg{\'e}briques
lisses sont gouvern{\'e}es par une {\'e}quation diff{\'e}rentielle,
dite de Picard-Fuchs. Pour celles-ci, et plus g{\'e}n{\'e}ralement,
pour leurs facteurs, Katz~\cite{katz72} puis Andr{\'e}~\cite{andre97}
relient les $p$-courbures {\`a} la r{\'e}duction modulo~$p$
de l'application de Kodaira-Spencer.
C'est ainsi qu'est {\'e}tablie la conjecture de Grothendieck
pour les {\'e}quations hyperg{\'e}om{\'e}triques de Gau{\ss} (on retrouve
la liste de Schwarz) ou les
facteurs des connexions de Knizhnik-Zamolodchikov.
L'approche arithm{\'e}tique est au c\oe ur de cet expos{\'e}.
Elle a permis
d'{\'e}tablir la conjecture pour les fibr{\'e}s de rang~$1$ sur une courbe
(Chudnovsky~\cite{chudnovsky85b},
le cas o{\`u} $E$ est trivial a {\'e}t{\'e} vu plus haut),
lorsque  $(E,\nabla)$ est extension de deux modules {\`a} connexion
isotriviaux (Andr{\'e}, \cite{andre97}),
et plus g{\'e}n{\'e}ralement, lorsque le \emph{groupe de Galois diff{\'e}rentiel}
de $(E,\nabla)$ a une composante neutre r{\'e}soluble (Andr{\'e}, voir plus bas).

Le groupe de Galois diff{\'e}rentiel {\og g{\'e}n{\'e}rique\fg}
d'un module {\`a} connexion $(E,\nabla)$
est d{\'e}fini par Katz dans~\cite{katz82} (voir aussi
l'expos{\'e}~\cite{bertrand92} de Bertrand {\`a} ce s{\'e}minaire).
Sa d{\'e}finition la plus {\'e}l{\'e}gante est tannakienne et fournit un groupe
alg{\'e}brique sur le corps $K(X)$ des fonctions de $X$: c'est le
sous-groupe de $\GL(E\otimes K(X))$  qui stabilise tout sous-module
horizontal de toute construction tensorielle sur $(E,\nabla)$.
Dans le cas singulier r{\'e}gulier, qui est en fin de compte celui
qui nous importe le plus pour la conjecture de Grothendieck,
on peut d{\'e}crire peut-{\^e}tre plus concr{\`e}tement une $\C$-forme de ce groupe.
{\'E}tant donn{\'e} un plongement de $K$ dans $\C$, on peut consid{\'e}rer
$(E,\nabla)$ comme un module {\`a} connexion int{\'e}grable sur la vari{\'e}t{\'e}
analytique $X(\C)$. Le groupe de Galois diff{\'e}rentiel {\og est\fg} alors
l'adh{\'e}rence de Zariski du groupe de monodromie usuel de $(E,\nabla)$
(\cite{katz82}, Prop.~5.2).
Dans tous les cas (\cite{katz82}, Prop.~4.5), 
la finitude du groupe de Galois diff{\'e}rentiel de $(E,\nabla)$
{\'e}quivaut {\`a} l'isotrivialit{\'e} de $(E,\nabla)$, soit encore
{\`a} l'existence d'une base de solutions alg{\'e}briques.

Katz a propos{\'e} (\cite{katz82}, Conj.~9.2) une g{\'e}n{\'e}ralisation
de la conjecture de Grothendieck et pr{\'e}dit que
l'alg{\`e}bre de Lie du groupe de groupe de Galois diff{\'e}rentiel
de $(E,\nabla)$ est la plus petite sous-alg{\`e}bre de Lie alg{\'e}brique
de $\mathfrak g\mathfrak l(E)$ qui {\og contient\fg}, pour presque
tout $p$, les $p$-courbures $\psi_p(\partial)$. 
En fait, et c'est l{\`a} l'objet principal de l'article~\cite{katz82},
cette g{\'e}n{\'e}ralisation et la conjecture de Grothendieck sont
{\'e}quivalentes: la v{\'e}racit{\'e} de la conjecture de Grothendieck
pour tout fibr{\'e} {\`a} connexion int{\'e}grable implique celle de Katz.

Enfin, signalons qu'un {\og $q$-analogue\fg}
de la conjecture de Grothendieck-Katz (pour les {\'e}quations
aux $q$-diff{\'e}rences)  a {\'e}t{\'e} d{\'e}montr{\'e} par L.~Di Vizio dans
sa th{\`e}se (d{\'e}cembre 2000).
La d{\'e}monstration utilise le th{\'e}or{\`e}me~\ref{theo.abdw}
ci-dessous, l'un des $R_v$ {\'e}tant infini.

\subsection{Exemples}
Revenons maintenant 
sur le th{\'e}or{\`e}me~\ref{theo.diff-exact} concernant l'exactitude de formes
diff{\'e}rentielles.
Soit $X$ une courbe alg{\'e}brique lisse sur un corps
de nombres~$K$ et soit $\omega\in\Gamma(X,\Omega^1_X)$ 
une forme diff{\'e}rentielle r{\'e}guli{\`e}re sur $X$. On peut consid{\'e}rer deux
modules {\`a} connexion (n{\'e}cessairement int{\'e}grables puisque $X$ est une courbe)
:
\begin{gather}
E_1 =\mathscr O_X, \qquad \nabla f  = \mathrm df - f \omega \\
E_2 = \mathscr O_X^2, \qquad \nabla (f,g) 
  =  (\mathrm df-\omega g, \mathrm dg). 
\end{gather}
L'existence d'une section horizontale non nulle de $E_1$ signifie
ainsi exactement que $\omega$ est une diff{\'e}rentielle logarithmiquement
exacte. Si $E_2$ est engendr{\'e} par ses sections horizontales,
il existe une telle section $(f,g)$ avec $g\neq 0$, donc $g$ constante
que l'on peut supposer {\'e}gale {\`a} $1$,
et alors $\mathrm df=\omega$, ce qui signifie que $\omega$
est une diff{\'e}rentielle exacte.
Dans les deux cas, l'hypoth{\`e}se du th{\'e}or{\`e}me~\ref{theo.diff-exact}
sur la r{\'e}duction modulo $p$ de $\omega$
pour presque tout nombre premier $p$ est, d'apr{\`e}s le
th{\'e}or{\`e}me de Cartier, {\'e}quivalente {\`a} la nullit{\'e} des $p$-courbures.

D'autre part, la d{\'e}finition m{\^e}me du groupe de Galois diff{\'e}rentiel
montre que c'est un sous-groupe de $\Gm$ dans le cas de
$(E_1,\nabla)$,
tandis que $(E_2,\nabla)$
ayant un sous-fibr{\'e} horizontal (d{\'e}fini par $q=0$),
son groupe de Galois diff{\'e}rentiel est un sous-groupe du groupe
de Borel $\left(\begin{smallmatrix} * & * \\ 0 &
*\end{smallmatrix}\right)$.
Ces deux groupes alg{\'e}briques sont en particulier r{\'e}solubles
et r{\'e}soudre la conjecture de Grothendieck pour les modules {\`a} connexion
int{\'e}grable dont le groupe de Galois est r{\'e}soluble d{\'e}montre
le th{\'e}or{\`e}me~\ref{theo.diff-exact}. R{\'e}ciproquement, si le groupe
de Galois  diff{\'e}rentiel est r{\'e}soluble, des d{\'e}vissages ram{\`e}nent
la conjecture de Grothendieck au cas du th{\'e}or{\`e}me~\ref{theo.diff-exact},
d'o{\`u} le th{\'e}or{\`e}me (\cite{andre97}, \S3,
voir aussi~\cite{bost2001}, th{\'e}or{\`e}me 2.9):

\begin{Theo}[Andr{\'e}]\label{theo.resoluble}
Soit $X$ une vari{\'e}t{\'e} alg{\'e}brique lisse sur un corps
de nombre et $(E,\nabla)$ un fibr{\'e} vectoriel
{\`a} connexion int{\'e}grable sur $X$.
Si la composante neutre du groupe de Galois diff{\'e}rentiel 
de $(E,\nabla)$ est r{\'e}soluble, alors $(E,\nabla)$ v{\'e}rifie
la conjecture de Grothendieck-Katz.
% \footnote{V{\'e}rifier pour Katz?}
% D'apres Katz, il suffit que la conj. de Grothendieck soit satisfaite pour
% toute construction de la forme constr(E) tensoriel le dual d'une
% construction de rang 1. Le gr de Galois diff d'un tel biniou est encore
% de composante neutre r{\'e}soluble.
\end{Theo}

\section{Feuilletages}\label{sec.feuilletages}

\subsection{}
Soit $X$ une vari{\'e}t{\'e} alg{\'e}brique lisse sur un corps $K$
et soit $F\subset T_X$ un sous-fibr{\'e} vectoriel
de son fibr{\'e} tangent qui est \emph{involutif}, c'est-{\`a}-dire
stable par le crochet de Lie sur $T_X$.
Si $K\subset\C$, un tel fibr{\'e} involutif d{\'e}finit un feuilletage
holomorphe sur la vari{\'e}t{\'e} analytique $X(\C)$: il existe dans un voisinage
$U_x$ de tout point $x\in X(\C)$ une sous-vari{\'e}t{\'e} analytique $Y\subset U_x$
telle que pour tout $y\in Y$, l'espace tangent $T_y Y$  {\`a} $Y$ en $y$
soit {\'e}gal {\`a} $F_y\subset T_y X$ (th{\'e}or{\`e}me de Frobenius).
Cependant, une telle \emph{feuille} n'est en g{\'e}n{\'e}ral pas un
germe de sous-vari{\'e}t{\'e} alg{\'e}brique.

Si $K$ est de caract{\'e}ristique~$0$ et $x\in X(K)$ un point rationnel,
la variante formelle du th{\'e}or{\`e}me de Frobenius fournit cependant une
sous-vari{\'e}t{\'e} formelle lisse de $X$ compl{\'e}t{\'e} le long de $x$,
la \emph{feuille formelle de $F$ en $x$}.

Si $K$ est de caract{\'e}ristique~$p$, 
on dira que $F$ est \emph{$p$-int{\'e}grable}
si c'est un sous-fibr{\'e} de~$\mathrm TX$ stable par puissance~\pieme.

Pour motiver cette d{\'e}finition,
supposons maintenant que $K$ soit un corps de nombres
et que la feuille formelle $\widehat Y$ passant par un point
rationnel~$x\in X(K)$ soit le germe d'une sous-vari{\'e}t{\'e} alg{\'e}brique~$Y$.
Pour presque tout nombre premier $p$,
on obtient par r{\'e}duction modulo~$p$ une situation analogue
$(X,F,Y)$ sur un corps fini de caract{\'e}ristique~$p$. 
Comme les op{\'e}rations de puissance
\pieme sur $\mathrm T_{Y}$ et sur $\mathrm T_X$ sont compatibles,
il en r{\'e}sulte que $F$ est clos par puissance~\pieme aux points
de~$Y$.

R{\'e}ciproquement, la g{\'e}n{\'e}ralisation de la conjecture de Grothendieck aux
feuilletages (formul{\'e}e par Ekedahl et Shepherd-Barron
dans~\cite{ekedahl-sb2000}, conjecture F)
pr{\'e}dit que les feuilles formelles
d'un tel sous-fibr{\'e} involutif $F\subset T_X$ qui pour presque
tout nombre premier $p$ est $p$-int{\'e}grable 
sont des sous-vari{\'e}t{\'e}s alg{\'e}briques.

Avant d'{\'e}noncer le th{\'e}or{\`e}me de Bost, quelques rappels 
d'analyse complexe s'imposent.

\subsection{}
Il est bien connu que toute fonction holomorphe born{\'e}e
sur l'espace affine $\C^n$ est constante: c'est le th{\'e}or{\`e}me
de Liouville.
Plus g{\'e}n{\'e}ralement, toute fonction \emph{plurisousharmonique}
major{\'e}e sur $\C^n$ est constante. Rappelons qu'une fonction
$\psi$ sur une vari{\'e}t{\'e} analytique complexe $M$ {\`a} valeurs
dans $\R\cup\{-\infty\}$
est dite plurisousharmonique~si: 
\begin{itemize}
\item elle est semicontinue sup{\'e}rieurement ;
\item elle n'est identiquement $-\infty$ sur aucune composante connexe de
$M$ ;
\item pour toute fonction holomorphe $f\colon D(0,1)\ra M$
du disque unit{\'e} ferm{\'e} de $\C$ dans $M$, on a l'in{\'e}galit{\'e}
\begin{equation}
  \psi(f(0))
 \leq \frac{1}{2\pi}
  \int _0^{2\pi} \psi(f(e^{\mathrm i\theta}))\,\mathrm d\theta.
\end{equation}
\end{itemize}
Par exemple, pour
toute fonction holomorphe $\phi$ non nulle, $\log\abs{\phi}$
est une fonction plurisousharmonique.

Suivant la terminologie de~\cite{bost2001},
nous dirons ainsi qu'une vari{\'e}t{\'e} analytique
complexe connexe $M$ v{\'e}rifie la \emph{propri{\'e}t{\'e} de Liouville}
si toute fonction plurisousharmonique major{\'e}e sur $M$
est constante.
Si $X$ est une vari{\'e}t{\'e} alg{\'e}brique complexe connexe lisse,
on peut d{\'e}montrer que $X(\C)$ satisfait la propri{\'e}t{\'e} de Liouville. 
Dans le cas compact, cela r{\'e}sulte du principe du maximum ;
dans le cas g{\'e}n{\'e}ral, on peut par exemple remarquer
que le compl{\'e}mentaire d'un ferm{\'e} analytique (strict) dans une
vari{\'e}t{\'e} qui satisfait la propri{\'e}t{\'e} de Liouville 
le satisfait aussi. Tout groupe de Lie complexe connexe
v{\'e}rifie la propri{\'e}t{\'e} de Liouville.

\begin{Theo}[Bost]\label{theo.feuilletage}
Soit $X$ une vari{\'e}t{\'e} alg{\'e}brique lisse sur un corps de nombres $K$
contenu dans $\C$,
soit $F$ un sous-fibr{\'e} involutif de $\mathrm TX$ et soit $x\in X(K)$ un point
rationnel.
Supposons que les deux conditions suivantes sont satisfaites:
\begin{itemize}
\item pour presque tout nombre premier $p$, la r{\'e}duction modulo~$p$
de $F$ est un sous-fibr{\'e} $p$-int{\'e}grable ;
\item 
il existe une vari{\'e}t{\'e} analytique complexe $M$ satisfaisant la propri{\'e}t{\'e}
de Liouville, un point~$O\in M$,
une application holomorphe $\psi $ de~$M$ vers
la feuille en~$x$ du feuilletage holomorphe induit par~$F$
sur $X(\C)$ telle que $\psi(O)=x$ et telle que $\psi$
soit biholomorphe d'un voisinage de $O$ vers un voisinage de $x$
dans cette feuille.
\end{itemize}
Alors, la feuille formelle de $F$ passant par $x$
est alg{\'e}brique.
\end{Theo}

\subsection{}
Les {\'e}quations diff{\'e}rentielles issues de fibr{\'e}s vectoriels {\`a} connexion
fournissent naturellement des feuilletages. Plus g{\'e}n{\'e}ralement,
si $G$ est un groupe alg{\'e}brique sur un corps de nombres $K$,
un $G$-fibr{\'e} principal {\`a} connexion sur une $K$-vari{\'e}t{\'e} alg{\'e}brique
lisse $B$ est un $G$-fibr{\'e} principal $X\ra B$
(c'est-{\`a}-dire un $G$-torseur sur $B$)
muni d'un scindage $G$-{\'e}quivariant de la suite exacte  de $G$-fibr{\'e}s
vectoriels sur $X$,
\begin{equation}
   0 \ra \mathrm T_{X/B} \ra \mathrm TX \ra \mathrm TB \ra 0.
\end{equation}
Dans ce cas, l'image de $\mathrm TB$ par la section fournit
un feuilletage $G$-{\'e}quivariant sur $X$.

{\`A} l'aide du th{\'e}or{\`e}me~\ref{theo.feuilletage},
on peut alors d{\'e}montrer certains cas de la conjecture de Grothendieck,
par exemple le th{\'e}or{\`e}me~\ref{theo.diff-exact} ou
le th{\'e}or{\`e}me~\ref{theo.resoluble}, si l'on {\'e}tablit
que la condition d'uniformisation est automatiquement
remplie.
Apr{\`e}s d{\'e}vissages, le th{\'e}or{\`e}me~\ref{theo.resoluble}
se ram{\`e}ne au cas o{\`u} la base $B$ est une courbe et
le groupe $G$ connexe commutatif.
Utilisant le morphisme d'addition $G^d\ra G$,
on peut alors {\'e}tendre le  $G$-fibr{\'e}  {\`a} connexion $(X,\nabla)$ sur $B$,
en un autre $(X_d,\nabla)$ sur la puissance
sym{\'e}trique $d$-i{\`e}me de $B$ (qui contient $B$
une fois fix{\'e} un point base).
Le point est que pour $d$ assez grand, le rev{\^e}tement universel
d'une telle puissance sym{\'e}trique v{\'e}rifie la propri{\'e}t{\'e} de Liouville.

\begin{Prop}
Soit $C$ une courbe alg{\'e}brique complexe lisse connexe.
Pour tout entier $d$ assez grand, le rev{\^e}tement universel
du produit sym{\'e}trique $\Sym^d C$ satisfait la propri{\'e}t{\'e} de Liouville.
\end{Prop}
Soit $\bar C$ la compl{\'e}tion projective lisse de~$C$,
$g$ son genre et $s$ le cardinal de $S=\bar C\setminus C$.
Il est d{\'e}montr{\'e} dans~\cite{bost2001} qu'il suffit
de prendre $d\geq \max(2,g,g+s-1)$. 
Cependant, si $d> 2g-2+s$, $\Sym^d C$ est un fibr{\'e} en espaces
affines (projectifs si $s=0$) au-dessus de la jacobienne
g{\'e}n{\'e}ralis{\'e}e $\Jac_S \bar C$ d{\'e}finie par le module~$S$.
Son rev{\^e}tement universel est l'image r{\'e}ciproque 
par l'application exponentielle $\Lie \Jac_S\bar C\ra \Jac_S\bar C$
de ce fibr{\'e} affine (resp.~projectif). Il v{\'e}rifie la condition
de Liouville. Si $S\neq\emptyset$, c'est m{\^e}me l'espace
affine de dimension~$d$.

\subsection{}
Dans~\cite{miyaoka87}, Miyaoka d{\'e}montre un th{\'e}or{\`e}me d'alg{\'e}bricit{\'e}
pour les feuilles de certains feuilletages alg{\'e}briques. 
Il proc{\`e}de essentiellement en \emph{construisant} des courbes
rationnelles tangentes au feuilletage.
Le point de d{\'e}part de la d{\'e}monstration est ainsi l'existence
d'une famille dense de courbes le long desquelles le fibr{\'e}
involutif d{\'e}finissant le feuilletage a un degr{\'e} strictement positif.
Les feuilles obtenues sont faiblement rationnellement connexes:
deux points quelconques peuvent {\^e}tre reli{\'e}s par une suite
de courbes rationnelles.

Le th{\'e}or{\`e}me~\ref{theo.feuilletage} admet un analogue g{\'e}om{\'e}trique,
d{\'e}couvert tr{\`e}s r{\'e}cemment de mani{\`e}re ind{\'e}pendante par
Bogomolov et McQuillan~\cite{bogomolov-mq2001} d'une part, 
en liaison avec les travaux de Miyaoka, et
par Bost, d'autre part.
L'{\'e}nonc{\'e} en est le suivant.
\begin{Theo}\label{theo.bogmq}
Soit $k$ un corps, disons alg{\'e}briquement clos
et soit $V_0$ une sous-vari{\'e}t{\'e} connexe projective lisse d'une
vari{\'e}t{\'e} quasi-projective~$X$ d{\'e}finie sur $k$.
On suppose que $V_0$ n'est pas r{\'e}duite {\`a} un point.
Soit $\widehat V$ une sous-vari{\'e}t{\'e} formelle lisse du compl{\'e}t{\'e} formel de $X$
le long de $V_0$ et contenant $V_0$.\footnote
{Lorsque $k=\C$ est le corps des nombres complexes,
une sous-vari{\'e}t{\'e} analytique
$V$ d'un voisinage ouvert $U$ de $V_0(\C)$ dans $X(\C)$ contenant
$V_0(\C)\cap U$ en fournit un excellent exemple.}
On suppose que le fibr{\'e} normal de $V_0$ dans $\widehat V$, $N_{V_0}\widehat V$
est ample. Alors, $\widehat V$ est alg{\'e}brique: l'adh{\'e}rence de Zariski
de $\widehat V$ dans $X$ a m{\^e}me dimension que $\widehat V$.
\end{Theo}
Autrement dit, il existe une sous-vari{\'e}t{\'e} alg{\'e}brique $Y$
de $X$ contenant $V_0$ dont $\widehat V$ est l'une des {\og branches\fg}
le long de~$V_0$.

La d{\'e}monstration est {\'e}tonnamment simple et, apr{\`e}s tout,
l'id{\'e}e n'est pas si diff{\'e}rente de celle qui 
conduira {\`a} la d{\'e}monstration du th{\'e}or{\`e}me~\ref{theo.feuilletage}.
Pr{\'e}cisons un peu les analogies entre ces deux th{\'e}or{\`e}mes.
Dans le th{\'e}or{\`e}me~\ref{theo.feuilletage},
il faut interpr{\'e}ter la vari{\'e}t{\'e} $X$
sur le corps de nombres~$K$ comme un germe de vari{\'e}t{\'e} $\mathscr X$ au-dessus
de la {\og courbe arithm{\'e}tique\fg} $\Spec\mathfrak o_K$. Le point
rationnel $x$ fournit alors une section, c'est-{\`a}-dire quelque chose
comme une courbe trac{\'e}e sur $\mathscr X$ qui est l'analogue
arithm{\'e}tique de la sous-vari{\'e}t{\'e} $V_0$ du th{\'e}or{\`e}me~\ref{theo.bogmq}.
Dans le th{\'e}or{\`e}me~\ref{theo.feuilletage}, l'alg{\'e}bricit{\'e} de la feuille
est plus simple {\`a} concevoir: dans ce cas, il existe effectivement
une vari{\'e}t{\'e} alg{\'e}brique dont le compl{\'e}t{\'e} formel en~$x$  s'identifie
{\`a} la feuille formelle passant par~$x$.

Rempla\c{c}ons $X$ par l'adh{\'e}rence de Zariski de $\widehat V$, c'est-{\`a}-dire
la plus petite sous-vari{\'e}t{\'e} alg{\'e}brique contenant $V_0$
et dont le compl{\'e}t{\'e} formel le long de $V_0$ contient $\widehat V$.
Notons $N=N_{V_0}\widehat V$,
$v_0=\dim V_0$ et $v=\dim\widehat V=v_0+\rg N$, de sorte que
la conclusion du th{\'e}or{\`e}me est que $v=\dim X$.
Soit $L$ la restriction {\`a} $X$ du fibr{\'e} en droites $\mathscr O(1)$
sur l'espace projectif. L'hypoth{\`e}se que $X$ est l'adh{\'e}rence
de Zariski de $\widehat V$ signifie que pour tout entier $D\geq 0$,
l'op{\'e}ration de restriction des sections de $L^D$ de $X$ {\`a} $\widehat V$
d{\'e}finit
une injection de $k$-espaces vectoriels
\begin{equation}\label{eq.injhatv}
\Gamma(X,L^D) \hra \Gamma(\widehat V,L^D).
\end{equation}
Nous allons filtrer l'espace vectoriel $E_D=\Gamma(\widehat V,L^D)$
par l'ordre d'annulation le long de $V_0$.
Si $n$ est un entier $\geq 1$, notons $V_n$ le 
voisinage infinit{\'e}simal d'ordre $n$ de $V_0$ dans $\widehat V$.
Notons $E_0^D=E^D$ et,  si $n\geq 1$,
soit $E_n^D$ l'ensemble des sections de $E_0^D$
dont la restriction {\`a} $V_{n-1}$ est nulle. Une section $s$ de $E_n^D$
poss{\`e}de un \emph{jet d'ordre~$n$} le long de $V_0$ ;
c'est une section sur $V_0$ de $L_{V_0}^D\otimes\Sym^n N^\vee$,
nulle si et seulement si $s$ appartient {\`a} $E_{n+1}^D$.
Par passage aux sous-quotients, on d{\'e}duit de l'injection~\eqref{eq.injhatv}
des injections, pour tous $n$ et $D$,
\begin{equation}
E_n^D/E_{n+1}^D  \hra \Gamma(V_0,L|_{V_0}^D \otimes \Sym^n N^\vee)
=  
\Gamma(\P(N^\vee), \pi^* L|_{V_0}^D\otimes\mathscr O_{\P(N^\vee)}(n)),
\end{equation}
o{\`u} l'on a not{\'e} $\pi\colon \P(N^\vee)\ra V_0$ le fibr{\'e} projectif associ{\'e} {\`a}
$N^\vee$ au-dessus de $V_0$.
Comme $\P(N^\vee)$ est une vari{\'e}t{\'e} projective sur $k$ de dimension
$v-1$, il existe une constante $c$ telle que pour tous
$D$ et $n$, on ait la majoration
\begin{equation}
\dim_k \Gamma(\P(N^\vee), \pi^* L^D\otimes\mathscr O_{\P(N^\vee)}(n))
  \leq c(1+n+D)^{v-1}.
\end{equation}     
D'autre part, l'hypoth{\`e}se d'amplitude sur $N$ et le th{\'e}or{\`e}me
d'annulation de Serre impliquent qu'il existe
un entier $\lambda\geq 1$ tel que si $n\geq \lambda D$
et $D\geq\lambda$,
\begin{equation}
\Gamma(V_0,L|_{V_0}^D \otimes \Sym^n N^\vee)
=0.
\end{equation}
Il en r{\'e}sulte que pour tout entier $D\geq\lambda$, on a l'in{\'e}galit{\'e}
\begin{equation}
\dim_k \Gamma(X,L^D) = \sum_{n=0}^\infty \dim_n E_n^D/E_{n+1}^D
 \leq \sum_{n=0}^{\lambda D-1} c_1(1+n+D)^{v-1}
       \leq c \lambda (1+\lambda)^{v-1} D^v .
\end{equation}
Par suite, $\dim_k\Gamma(X,L^D)=\mathrm O(D^v)$. D'apr{\`e}s
le th{\'e}or{\`e}me de Hilbert, on a $\dim X=v$,
ainsi qu'il fallait d{\'e}montrer.

\section{Fragments de th{\'e}orie d'Arakelov}\label{sec.arakelov}

Nous r{\'e}sumons maintenant
les quelques d{\'e}finitions et r{\'e}sultats de th{\'e}orie d'Arakelov
dont nous aurons besoin. 
Le lecteur int{\'e}ress{\'e} trouvera des compl{\'e}ments importants
dans les r{\'e}f{\'e}rences~\cite{szpiro85,soule-a-b-k92,bost96, soule97}.

\subsection{}
Soit $K$ un corps de nombres ; notons $\mathfrak o_K$
son anneau d'entiers. Les valeurs absolues sur $K$ sont de deux types.
Celles qui sont non-archim{\'e}diennes (on dit aussi finies)
sont associ{\'e}es {\`a}
un id{\'e}al premier $\mathfrak p$ de $\mathfrak o_K$ ;
on note $v_{\mathfrak p}:K^\times\ra\Z$,
la valuation associ{\'e}e et $\abs{\cdot}_{\mathfrak p}$
la valeur absolue correspondante, normalis{\'e}es de sorte que,
si $\pi$ est une uniformisante, $v_{\mathfrak p}(\pi)=1$
et $\log \abs{\pi}_{\mathfrak p}=- \log\#(\mathfrak o_K/\mathfrak p)$.
Si $K_{\mathfrak p}$ d{\'e}signe le compl{\'e}t{\'e} $\mathfrak p$-adique
de $K$ et $p$ la caract{\'e}ristique du corps r{\'e}siduel $\mathfrak
o_K/\mathfrak p$, on a ainsi
\begin{equation}
 \log \abs{p}_{\mathfrak p} = - [K_{\mathfrak p}:\Q_p] \log p. 
\end{equation}
Les valeurs absolues archim{\'e}diennes (ou infinies)
sont associ{\'e}es {\`a} un plongement\break $\sigma\colon K\hra\C$,
mais deux plongements conjugu{\'e}s fournissent la m{\^e}me valeur absolue.
Si $\sigma\colon K\hra \C$ est un tel plongement, on 
d{\'e}finit $\eps_{\sigma}=1$ si $\sigma$ est r{\'e}el et $2$ sinon,
et on a, si $x\in K$,
\begin{equation}
  \log\abs{x}_\sigma = \eps_\sigma \log\abs{\sigma(x)}. 
\end{equation}
On notera aussi $K_\sigma=\R$ ou $\C$ suivant que $\sigma$
est r{\'e}el ou complexe. On a ainsi $\eps_\sigma=[K_\sigma:\R]$.

Notons $\Sigma_K$ l'ensemble des valuations sur $K$ ainsi d{\'e}finies;
on notera aussi $\Sigma_{K,f}$ et $\Sigma_{K,\infty}$ les ensembles
de valuations non-archim{\'e}diennes et archim{\'e}diennes respectivement.
Avec ces normalisations, on a la \emph{formule du produit}:
\begin{equation}
\text{si $x\in K\setminus\{0\}$,}\qquad
 \prod_{v\in\Sigma_K} \abs{x}_v = 1.
\end{equation}

\begin{defi}
Un $\mathfrak o_K$-fibr{\'e} vectoriel hermitien $\bar E=(E,\norm{\sigma})$
est la donn{\'e}e d'un $\mathfrak o_K$-module
projectif de rang fini $E$, ainsi que pour tout $\sigma\in
\Sigma_{K,\infty}$,
d'une norme hermitienne 
sur l'espace vectoriel $E_\sigma=\C\otimes_{\mathfrak o_K,\sigma} E$,
invariante par la conjugaison complexe si $\sigma$ est r{\'e}elle.
\end{defi}
Il revient au m{\^e}me de consid{\'e}rer une norme hermitienne
invariante par conjugaison 
sur l'espace vectoriel complexe $\C\otimes_{\Z} E$.

Soit $\bar E$ un $\mathfrak o_K$-fibr{\'e} vectoriel hermitien.
Si $\mathfrak p$ est un id{\'e}al premier de $\mathfrak o_K$,
on peut aussi d{\'e}finir une norme $\mathfrak p$-adique sur
l'espace vectoriel $E_{\mathfrak p}=K_{\mathfrak p}\otimes_{\mathfrak o_K} E$:
par d{\'e}finition, si $e\in E_\mathfrak p$,
\begin{equation}
\norm{e}_{\mathfrak p} = \inf \{ \abs{a}_\mathfrak p\,;\,
             a\in K_\mathfrak p^\times,\quad ae\in
\mathfrak o_{\mathfrak p} 
\otimes_{\mathfrak o_K} E \}.
\end{equation}

\subsection{}
Plus g{\'e}n{\'e}ralement, si $\mathscr X$ est un $\Z$-sch{\'e}ma
de type fini et plat, un \emph{fibr{\'e} vectoriel hermitien}
$\bar E=(E,\norm{})$
sur $\mathscr X$ est la donn{\'e}e d'un fibr{\'e} vectoriel $E$ sur $\mathscr X$
ainsi que d'une \emph{m{\'e}trique hermitienne} 
sur le fibr{\'e} vectoriel holomorphe
$E_\C$ sur l'espace analytique complexe  $\mathscr X(\C)$,
suppos{\'e}e invariante par conjugaison complexe.

On peut effectuer, sur ces fibr{\'e}s vectoriels hermitiens, un certain
nombre de constructions standard: image r{\'e}ciproque
par un morphisme $\mathscr X'\ra\mathscr X$; somme directe,
munie de la norme somme directe orthogonale;
sous-fibr{\'e}, muni de la norme induite; fibr{\'e} quotient,
muni de la norme quotient; modules d'homomorphismes, en particulier duaux;
produit tensoriel, puissances sym{\'e}triques et ext{\'e}rieures (toutes
deux \emph{quotients} du produit tensoriel).

\subsection{}\label{arak.globalsect}
Soit $\mathscr X$ un $\mathfrak o_K$-sch{\'e}ma de type fini
et plat.
Si $\sigma$ est un plongement de $K$ dans $\C$, notons
$\mathscr X_\sigma(\C)=\mathscr X\otimes_{\mathfrak o_K,\sigma} \C$,
de sorte que $\mathscr X(\C)$ est la r{\'e}union disjointe des
$\mathscr X_{\sigma}(\C)$.
Si $\bar E$ est un fibr{\'e} vectoriel hermitien sur $\mathscr X$
et si $\mathscr X$ est propre, 
$\Gamma(\mathscr X,E)$ est un $\mathfrak o_K$-module projectif
de rang fini. On peut le munir d'une structure de $\mathfrak o_K$-fibr{\'e}
vectoriel hermitien comme suit. Fixons $\mu$ une mesure {\og de
Lebesgue\fg} sur $\mathscr X(\C)$ (voir~\cite{bost2001}
pour la d{\'e}finition pr{\'e}cise;
quand $\mathscr X(\C)$ est lisse
de dimension complexe~$d$
c'est une mesure qui, en coordonn{\'e}es locales $z_j=x_j+\mathrm iy_j$,
s'{\'e}crit
$\omega(x) \mathrm dx_1\mathrm dy_1\dots\mathrm dx_d\mathrm dy_d$
o{\`u} $\omega$ est une fonction $\mathscr C^\infty$ strictement positive).
Alors, on d{\'e}finit pour tout $e\in \Gamma(\mathscr X,E)$
et tout plongement complexe $\sigma\colon K\hra\C$,
\begin{equation}
   \norm{e}_\sigma = \left(  \int_{\mathscr X_\sigma(\C)}
        \norm{e(x)}_\sigma^2\, \mathrm d\mu(x) \right)^{1/2}.
\end{equation}
On d{\'e}finit aussi une \emph{norme sup.}:
\begin{equation}
\norm{e}_{\infty,\sigma}
       = \sup_{x\in\mathscr X_\sigma(\C)} \norm{e(x)}_{\sigma}.
\end{equation}

Un lemme {\'e}l{\'e}mentaire relie ces deux normes.
\begin{Lemm}\label{lemm.gromov}
Si de plus $\bar E$ est un fibr{\'e} en droites hermitien sur $\mathscr X$,
il existe une constante $C$ telle que pour tout $\sigma\colon K\hra\C$,
tout entier $D\geq 1$ et tout $e\in\Gamma(\mathscr X,E^D)$,
on ait l'encadrement
\begin{equation} \label{ineq.gromov}
 \norm{e}_\sigma\leq C \norm{e}_{\infty,\sigma}  \leq C^D \norm{e}_\sigma .
\end{equation}
\end{Lemm}

\subsection{}
Un $\mathfrak o_K$-fibr{\'e} vectoriel hermitien $\bar E$ de rang~$1$
poss{\`e}de un \emph{degr{\'e}}, d{\'e}fini par 
\begin{equation}
 \hdeg \bar E = \log \# (E/\mathfrak o_K e)
   - \sum_{\sigma\colon K\hra\C} \log \norm{e}_\sigma ,
\end{equation}
$e$ {\'e}tant un {\'e}l{\'e}ment non nul quelconque de $E$;
le fait que la formule n'en d{\'e}pende pas r{\'e}sulte de la formule du produit.
Plus g{\'e}n{\'e}ralement, pour tout {\'e}l{\'e}ment non nul $e\in  E_K$,
on a
\begin{equation}
 \hdeg\bar E = - \sum_{v\in \Sigma_K} \log \norm{e}_v .
\end{equation}
On d{\'e}finit plus g{\'e}n{\'e}ralement le degr{\'e} d'un $\mathfrak o_K$-fibr{\'e} vectoriel
hermitien $\bar E$ de rang $d\geq 1$ comme celui de sa puissance ext{\'e}rieure
maximale $\bigwedge^d \bar E$.
On d{\'e}finit aussi la \emph{pente} de $\bar E$, not{\'e}e $\hmu(\bar E)$, 
comme le quotient de son degr{\'e} par son rang.

Si $\bar E$, $\bar E'$ et $\bar F$ sont trois $\mathfrak o_K$-fibr{\'e}s vectoriels
hermitiens, $\bar F$ {\'e}tant un sous-fibr{\'e} de $\bar E$, on a les formules
\begin{gather}
\hdeg \bar E^\vee = -\hdeg \bar E , \\
\label{eq.degadditif}
\hdeg \bar E = \hdeg \bar F + \hdeg (\bar E/\bar F), \\
\hmu (\bar E\otimes \bar E') = \hmu (\bar E)+ \hmu(\bar E').
\end{gather}

D'apr{\`e}s le th{\'e}or{\`e}me des facteurs invariants et la th{\'e}orie
des d{\'e}terminants de Gram,
si $(e_1,\dots,e_d)$ est une famille lin{\'e}airement
ind{\'e}pendante d'{\'e}l{\'e}ments de $E$, on a
\begin{equation}
\hdeg \bar E = \log \#(E/ (\mathfrak o_K e_1+\dots+\mathfrak o_K e_d))
      - \frac12 \sum_{\sigma\colon K\hra\C} 
               \log \det \big( \langle e_j,e_k\rangle_\sigma \big).
\end{equation}
En utilisant l'in{\'e}galit{\'e} d'Hadamard, on en d{\'e}duit que pour
toute famille g{\'e}n{\'e}ratrice $(e_1,\dots,e_r)$ de $K\otimes E$,
\begin{equation}\label{eq.hadamard2}
\hmu(\bar E)
 \geq - \sum_{v\in\Sigma_K} \log \max_{1\leq j\leq r} \norm{e_j}_v.
\end{equation}
 
\begin{Prop}\label{prop.hsa}
Soit $\mathscr X$ un $\mathfrak o_K$-sch{\'e}ma propre et plat
et soit $\bar L$ un fibr{\'e} en droites hermitien sur $\mathscr X$.
Si $D$ est un entier $\geq 0$, notons $\bar E_D$
le $\mathfrak o_K$-fibr{\'e} vectoriel hermitien
d{\'e}fini par $\Gamma(\mathscr X,L^D)$.
Si $L$ est ample, il existe une constante $c\in\R$
telle que pour tout entier $D\geq 0$, on ait
\begin{equation}
\hmu (\bar E_D) \geq - c D.
\end{equation}
\end{Prop}
Pour la d{\'e}monstration, combiner l'in{\'e}galit{\'e}~\eqref{eq.hadamard2},
le lemme~\ref{lemm.gromov} et le fait que la\break $K$-alg{\`e}bre gradu{\'e}e
$\bigoplus_D \Gamma(\mathscr X_K,L^D)$  est de type fini.

En appliquant l'in{\'e}galit{\'e}~\eqref{eq.hadamard2} au fibr{\'e}
dual $\bar E^\vee$, on d{\'e}montre
que lorsque $\bar F$
parcourt les sous-fibr{\'e}s vectoriels hermitiens de $\bar E$,
l'ensemble des pentes $\hmu(\bar F)$ est major{\'e}. La borne
sup{\'e}rieure (en fait, un maximum) est la \emph{plus grande pente}
de $\bar E$ et est not{\'e}e $\hmumax(\bar E)$.
Le comportement pr{\'e}cis de $\hmumax$ par produit tensoriel
est d{\'e}licat.
Si $\bar E$ et $\bar L$ sont deux $\mathfrak o_K$-fibr{\'e}s
vectoriels hermitiens, $\bar L$ {\'e}tant de rang~$1$, on a
\begin{equation}\label{eq.hmumax(etensl)}
\hmumax(\bar E\otimes\bar L) =\hmumax(\bar E) + \hdeg \bar L
\end{equation}
Pour des puissances sym{\'e}triques,
on {\'e}tablit 
le comportement asymptotique:
\begin{Lemm}\label{lemm.symk}
Soit $\bar E$ un $\mathfrak o_K$-fibr{\'e} vectoriel hermitien.
Il existe une constante $c\in\R$ telle que pour tout entier $k\geq 0$,
\begin{equation}
 \hmumax(\Sym^k \bar E) \leq c k.
\end{equation}
\end{Lemm}
Soit $\bar F$ un sous-fibr{\'e} de $\Sym^k \bar E$.
Il correspond par dualit{\'e} {\`a} un quotient $\bar G$ du
fibr{\'e} hermitien $(\Sym^k \bar E)^\vee$,
lequel s'identifie au sous-fibr{\'e} hermitien $\Gamma^k\bar E^\vee$
de $(\bar E^\vee)^{\otimes k}$ 
fix{\'e} par l'action du groupe sym{\'e}trique $\mathfrak S_k$.
Soit $(e_1,\dots,e_r)$ une base de $E^\vee_K$.
Si $\mathbf n=(n_1,\dots,n_r)\in\N^r$ avec $\sum_{i=1}^r n_i=k$,
soit $f_{\mathbf n}$ l'orbite
du tenseur $e_1^{n_1}\otimes\dots\otimes
e_r^{n_r}$ sous l'action de $\mathfrak S_k$.
Pour toute place finie $v$ de $K$, on a ainsi
\[ \norm{f_{\mathbf n}}_v \leq \max (\norm{e_1}_v,\dots,\norm{e_r}_v)^k
\]
tandis que pour toute place archim{\'e}dienne~$v$,
\[ \norm{f_{\mathbf n}}_v \leq \frac{k!}{n_1!\dots n_r!}
        \max (\norm{e_1}_v,\dots,\norm{e_r}_v)^k
 \leq r^k \max (\norm{e_1}_v,\dots,\norm{e_r}_v)^k.
\]
Les $f_{\mathbf n}$ forment une base de $\Gamma^k E^\vee_K$.
L'in{\'e}galit{\'e}~\eqref{eq.hadamard2} appliqu{\'e}e {\`a} leurs images
dans $\bar G$ jointe {\`a} l'{\'e}galit{\'e} $\hmu(\bar G)=-\hmu(\bar F)$
implique ainsi que
\[ \hmu (\bar F) \leq  k \left(
  [K:\Q]\log r + 
\sum_{v\in\Sigma_K}
\log\max(\norm{e_1}_v,\dots,\norm{e_r}_v) 
\right), \]
comme il fallait d{\'e}montrer.

\begin{Rema}
En utilisant les r{\'e}sultats de Zhang~\cite{zhang95}
qui fournissent une {\og presque-base orthonorm{\'e}e\fg}
de $\bar E$,
Bost a donn{\'e} une version explicite du lemme~\ref{lemm.symk},
cf.~\cite{graftieaux2001}.
Par ailleurs, dans la situation de la proposition~\ref{prop.hsa},
l'existence et le calcul de la limite
de $\hmu(\bar E_D)/D$, lorsque $D\ra +\infty$,
est le~\emph{th{\'e}or{\`e}me de Hilbert-Samuel arithm{\'e}tique},
cf.~\cite{gillet-s88,abbes-b95,rumely-l-v2000}.
\end{Rema}

\subsection{}
Soit $\bar E$ et $\bar F$ deux $\mathfrak o_K$-fibr{\'e}s vectoriels hermitiens
et soit $\phi\colon E_K\ra F_K$ une application $K$-lin{\'e}aire
injective.
Pour toute valuation $v\in\Sigma_K$, d{\'e}finissons la \emph{hauteur de
$\phi$ en $v$} par
\begin{equation}
h_v(\phi) = \log \norm{\phi}_v
= \log \sup_{e\in E_v\setminus\{0\}} \frac{\norm{\phi(e)}_v}{\norm{e}_v}.
\end{equation}
Pour tout $v$ sauf pour un nombre fini, $h_v(\phi)=0$, si bien
qu'on peut d{\'e}finir
la \emph{hauteur} (globale) de $\phi$ par  la formule
\begin{equation}
h(\phi) = \sum_{v\in\Sigma_K} h_v(\phi).
\end{equation}
En combinant les d{\'e}finitions et l'in{\'e}galit{\'e} de Hadamard,
on {\'e}tablit alors \emph{l'in{\'e}galit{\'e} de pentes} sur laquelle toute l'histoire
qui va suivre est fond{\'e}e:
\begin{equation}
\hmu( \bar E) \leq \hmumax(\bar F) + h(\phi).
\end{equation}
Plus pr{\'e}cis{\'e}ment, nous aurons besoin de la variante {\og filtr{\'e}e\fg}
de cette in{\'e}galit{\'e} de pentes, dont la d{\'e}monstration
est imm{\'e}diate {\`a} partir de l'in{\'e}galit{\'e} pr{\'e}c{\'e}dente et de la
formule~\eqref{eq.degadditif}.
\begin{Prop}\label{prop.pentes}
Soit $\bar E$ et $(\bar G^{(n)})_{n\geq 0}$ des
$\mathfrak o_K$-fibr{\'e}s vectoriels hermitiens.
Soit $ F_K$ un $K$-espace vectoriel
muni d'une filtration
d{\'e}croissante s{\'e}par{\'e}e $(F^{(n)}_K)_{n\geq 0}$
et, pour tout entier $n\geq 0$,
un isomorphisme $F^{(n)}_K/F^{(n+1)}_K\simeq G^{(n)}_K$.

Soit $\phi\colon E_K\ra F_K$ une application lin{\'e}aire injective.
Pour tout $n\geq 0$,
soit $\bar E^{(n)}$ le\break $\mathfrak o_K$-fibr{\'e} vectoriel
hermitien d{\'e}fini par $E^{(n)}=\phi^{-1}(F^{(n)}_K)$, muni des normes
induites par $\bar E$ et soit $\phi^{(n)}_K$ l'application
lin{\'e}aire induite $E^{(n)}_K\ra G^{(n)}_K$.

Alors, on a l'in{\'e}galit{\'e} 
\begin{equation}
 \hdeg\bar E \leq \sum_{n=0}^\infty \rg (E^{(n)}/E^{(n+1)}) \left(
             \hmumax(\bar G^{(n)}) + h(\phi^{(n)}) \right). 
\end{equation}
\upshape{(C'est une somme finie si l'on convient que
l'expression entre parenth{\`e}ses est nulle lorsque que $E^{(n)}=E^{(n+1)}$.)}
\end{Prop}

\section{Un premier th{\'e}or{\`e}me d'alg{\'e}bricit{\'e}}\label{sec.abdw}

L'article d'Andr{\'e}~\cite{andre97} (voir aussi~\cite{andre89})
repose sur un th{\'e}or{\`e}me d'alg{\'e}bricit{\'e} d'une s{\'e}rie formelle,
sorte d'analogue du th{\'e}or{\`e}me de Borel-Dwork.
Nous donnons ici une d{\'e}monstration de ce th{\'e}or{\`e}me
qui utilise le formalisme introduit au chapitre pr{\'e}c{\'e}dent.
Pour all{\'e}ger les notations, nous nous contentons d'une variable
quoique la g{\'e}n{\'e}ralisation du crit{\`e}re de P\'olya cit{\'e} {\`a}
la fin du chapitre~\ref{sec.4thms} en n{\'e}cessite plusieurs.

\subsection{}
Soit $K$ un corps de nombres et soit $y=\sum_{n=0}^\infty a_n x^n$ une s{\'e}rie
formelle {\`a} coefficients dans $K$.
On d{\'e}finit trois invariants, $\rho$, $\sigma$ et $\tau$:
\begin{gather}
\rho(y) = \sum_{v\in\Sigma_K}
     \limsup_{n\ra\infty} \frac1n\sup_{m\leq n}
                \log^+\abs{a_m}_v,  \\
\sigma(y) = \limsup_{n\ra\infty} \sum_{v\in\Sigma_K} \frac1n \sup_{m\leq n}
         \log^+\abs{a_m}_v , \\
\tau(y) = \inf_{\substack{S\subset\Sigma_K \\ \text{$S$ finie}}}
                 \limsup_{n\ra\infty} \sum_{v\in\Sigma_K\setminus S}
\frac1n \sup_{m\leq n} \log^+\abs{a_m}_v.
\end{gather}
(Si $x\in\R$, $\log^+(x)=\log \max(1,x)$.)
Si $S$ est une partie de $\Sigma_K$, $\rho_S(y)$ et $\sigma_S(y)$ d\'esigneront
les quantit\'es analogues \`a $\rho(y)$ et  $\sigma(y)$ o\`u 
seules les places de~$S$ sont prises en compte. 
Par exemple, il vient ainsi $\tau(y)=\inf_S \sigma_S(y)$.
Notons de plus $R_v(y)$ le rayon de convergence $v$-adique
de la s\'erie $y$ ; on a pour toute ensemble
$S$ de places l'\'egalit\'e $\rho_S(y)=\sum_{v\in S} \log^+ R_v(y)^{-1}$.

S'il existe un ensemble fini de places $S\subset\Sigma_K$
telles que les coefficients de $y$ soient $S$-entiers,
on a $\tau(y)=0$. Dans ce cas, $\rho(y)<\infty$ si et seulement
si pour toute place $v\in S$, le rayon de convergence $v$-adique
de la s\'erie $y$ n'est pas nul.

Soit $v$ une place de $K$. Une \emph{uniformisation
$v$-adique simultan{\'e}e} de $x$ et $y$ dans le disque  {\og ouvert\fg}
$D(0,R_v)\subset\bar K_v$
est la donn{\'e}e de deux fonctions m{\'e}romorphes $v$-adiques
$\phi$ et $\psi$ v{\'e}rifiant
\begin{itemize}
\item $\phi(0)=0$ et $\phi'(0)=1$ ;
\item $y(\phi(z))$ est le germe en l'origine de la fonction m{\'e}romorphe
 $\psi(z)$.
\end{itemize}
C'est ainsi, en quelque sorte, une uniformisation (m{\'e}romorphe)
du graphe de $y$ par le disque $D(0,R_v)$.
(Par d{\'e}finition, une fonction m{\'e}romorphe sur $D(0,R_v)$ est
le quotient de deux fonctions analytiques sur ce disque.)

On parle d'\emph{uniformisation triviale}
si de plus $\phi$ est l'identit{\'e}.

La premi{\`e}re partie du th{\'e}or{\`e}me suivant est due {\`a} Andr{\'e}
(\cite{andre97}, th{\'e}or{\`e}me 2.3.1), la seconde est essentiellement
le crit{\`e}re de Borel-Dwork.
\begin{Theo}\label{theo.abdw}
Soit $y\in K[[x]]$ telle que $\tau(y)=0$ et $\rho(y)<\infty$.
On suppose
que pour toute place $v$ de $K$, il existe une uniformisation
$v$-adique simultan{\'e}e de $x$ et $y$ dans un disque
$D(0,R_v)$. Si $\prod R_v>1$, alors $y$ est une fonction alg{\'e}brique.

Si de plus les uniformisations sont triviales pour tout $v$,
alors $y$ est une fonction rationnelle.
\end{Theo}

La fin de ce chapitre est consacr{\'e}e {\`a} la d{\'e}monstration du
th{\'e}or{\`e}me~\ref{theo.abdw}.

\subsection{}
Soit $d$ et $D$ deux entiers $\geq 1$
et soit $E_{d,D}\subset \mathfrak o_K[X,Y]$ le $\mathfrak o_K$-module
libre des polyn{\^o}mes de degr{\'e}s $\leq d$ en $X$ et $\leq D$ en $Y$.
On le munit des normes hermitiennes induites par la
base standard aux places archim{\'e}diennes, d'o{\`u} un 
$\mathfrak o_K$-fibr{\'e} vectoriel hermitien $\bar E_{d,D}$
de rang $(d+1)(D+1)$ et de degr{\'e} arithm{\'e}tique nul.

Soit $F_K=K[[x]]$ et soit $\phi\colon E_{d,D;K}\ra F_K$
l'application lin{\'e}aire d{\'e}finie par $P\mapsto P(x,y(x))$.
En filtrant $K[[x]]$ par l'ordre d'annulation en l'origine,
soit $F_K^{(k)}=x^k K[[x]]$,  on est dans la situation de la
proposition~\ref{prop.pentes}, o{\`u} pour tout $k\geq 0$,
$\bar G^{(k)}= \mathfrak o_K$ muni de la norme triviale
$\norm{1}_\sigma=1$ ; en particulier, $\hdeg \bar G^{(k)}=
\hmumax( \bar G^{(k)})=0$.

On raisonne par l'absurde.
Supposons que $y$ n'est pas une fonction alg{\'e}brique.
Alors, pour tous $d$ et $D$, l'application lin{\'e}aire $\phi$
est injective.
Si $y$ n'est pas une fonction rationnelle, l'application
$\phi$ est injective pour $D=1$ et tout entier $d\geq 1$.
L'in{\'e}galit{\'e} de pentes de la proposition~\ref{prop.pentes} 
s'{\'e}crit ainsi
\begin{equation}\label{eq.inegpentes}
0 \leq \sum_{n=0}^\infty \rg (E_{d,D}^{(n)}/E_{d,D}^{(n+1)})
           h(\phi^{(n)}). 
\end{equation}
Le reste de la d{\'e}monstration consiste {\`a} majorer convenablement
$h(\phi^{(n)})$ de sorte {\`a} contredire 
l'in{\'e}galit{\'e} pr{\'e}c{\'e}dente.
Pour cela, on utilise deux types d'estimations: 1\textsuperscript o)~pour
presque toute place, une majoration
{\og triviale\fg} (lemme~\ref{Lemm.evidente})
qui repose sur les hypoth\`eses $\tau(y)=0$ et
$\rho(y)<\infty$, et 2\textsuperscript o)~pour un nombre
fini de places, une majoration fond\'ee sur le lemme
de Schwarz, sous la forme:
\begin{Lemm}
Soit $v$ une place de $K$ et soit $f$ une fonction analytique
born{\'e}e sur le disque $D(0,R_v)\subset \bar K_v$.
Si $f$ s'annule {\`a} l'ordre $n$ en $0$, on a
\[ \abs{\frac{1}{n!} f^{(n)}(0) }_v
           \leq R_v^{-n} \norm{f}_{R_v}.
\]
\end{Lemm}

\subsection{}
Pla\c{c}ons-nous dans le cas d'une uniformisation triviale en une place
$v$ de $K$ et
supposons qu'il existe une fonction m{\'e}romorphe $\phi$
sur le disque $D(0,R_v)$ dont $y$ soit le d{\'e}veloppement de Taylor
en l'origine. Si $R'_v<R_v$, il existe alors des fonctions
analytiques born{\'e}es sur le disque $D(0,R'_v)$ telles que 
$\phi=f/g$ et $g(0)=1$.
Soit $n\geq 0$ et soit $P\in E_{d,D}^{(n)}$.
Si $P=\sum_{i=0}^D P_i(X) Y^i$, o{\`u} pour tout $i$, $\deg P_i\leq d$,
d{\'e}finissons une fonction analytique born{\'e}e $h$ sur $ D(0,R'_v)$ par
\[ h(x) = g(x)^D P(x,y(x)) = \sum_{i=0}^D P_i(x) f(x)^i g(x)^{D-i}.
\]
Comme $g(0)=1$ et comme $P(x,y(x))$ est suppos{\'e} {\^e}tre d'ordre au moins
$n$ en l'origine, on a
\[ \phi^{(n)}(P)
 = \frac{1}{n!} \frac{\partial^n}{\partial x^n} P(x,y(x))|_{x=0}
 = \frac{1}{n!} \frac{\partial^n}{\partial x^n} h(x)|_{x=0}.
\]
D'apr{\`e}s le lemme de Schwarz, on a alors
\[ \abs{\phi^{(n)}(P)}_v
 \leq R_v^ {\prime-n} R_v^{\prime d}
      \max(\norm{f}_{R'_v},\norm{g}_{R'_v})^D
	  \norm{P}_v
\]
multipli{\'e} par $\rg E_{d,D}$ si $v$ est une place archim{\'e}dienne.
On en d{\'e}duit qu'il existe une constante $C_v$ telle que l'on ait
pour tous $n$, $d$, $D$, l'in{\'e}galit{\'e}
\begin{equation}\label{eq.schwarz-triv}
 h_v(\phi^{(n)} )\leq (d-n) \log R'_v + D C_v
\end{equation}
auquel il faut ajouter $\log\rg E_{d,D}$ si $v$ est archim{\'e}dienne.

\subsection{}
Dans le cas d'une uniformisation simultan{\'e}e g{\'e}n{\'e}rale, le m{\^e}me
raisonnement fournit l'existence pour tout r{\'e}el $R'_v<R_v$
d'une constante $C_v$
telle que l'on ait pour tous $n$, $d$, $D$, l'in{\'e}galit{\'e}
\begin{equation}\label{eq.schwarz-simul}
 h_v(\phi^{(n)}) \leq -n \log R'_v + (d+D) C_v
\end{equation}
auquel il faut encore ajouter $\log\rg E_{d,D}$ si $v$ est archim{\'e}dienne.

\begin{Lemm} \label{Lemm.evidente}
Pour tout ensemble de places $T\subset\Sigma_K$, et pour
tout $\eps>0$, il existe une constante $C_T(\eps)$
telle que l'on ait la majoration 
\begin{multline*} \sum_{v\in T} h_v(\phi^{(n)})
       \leq  n \big(\rho_T(y)+\eps(1+\log D)\big)
      + C_T(\eps)  \\ 
   + [K:\Q]\log\rg E_{d,D} + [K:\Q](D+1)\log n .
\end{multline*}
(On peut omettre les deux derniers termes si $T$ ne contient
pas de places archim\'ediennes.)
\end{Lemm}
Si $j\geq 0$,
introduisons le d{\'e}veloppement en s{\'e}rie de $y^j$,
$y=\sum_{m=0}^\infty a_m^{(j)} x^m$, o{\`u} les $a_m^{(j)}$ sont dans $K$.
Constatons que l'application $\phi^{(n)}$
associe {\`a} un polyn{\^o}me $P\in E_{d,D;K}$ tel que $P(x,y(x))\in x^nK[[x]]$
le coefficient de $x^n$ dans $P(x,y(x))$.
Ainsi, si $v$ est une place finie, on a l'estim{\'e}e {\'e}vidente:
\begin{equation}% \label{eq.evidente}
        h_v(\phi^{(n)}) \leq
   \sup_{\substack { m\leq n \\ j\leq D }}
     \log^+\abs{a_m^{(j)}}_v.
\end{equation}
Lorsque $v$ est une place archim\'edienne, il faut rajouter 
$ \log\rg E_{d,D} + (D+1)\log n$.

Par suite, pour toute place $v$ de~$K$,
\[ \limsup \frac1n h_v(\phi^{(n)})
      \leq \sup_{j\leq D}
      \limsup_n \sup_{m\leq n} \frac1n \log^+\abs{a_m^{(j)}}_v
   \leq \sup_{j\leq D} \log^+ R_v(y^j)^{-1}, \]
$R_v(y^j)$ d\'esignant le rayon de convergence $v$-adique de la s\'erie $y^j$.
Mais ce rayon n'est autre que celui de $y$, si bien que
\[ \limsup \frac1n h_v(\phi^{(n)}) \leq \log^+ R_v(y)^{-1}. \]

D'autre part, si $v$ est une place finie de~$K$ et si $j\leq D$,
$a_m^{(j)}$ est major\'e par le maximum des produits
$a_{m_1}\dots a_{m_j}$, pour $m_1+\dots+m_j=j$. Si ce maximum
est atteint en $(m_1,\dots,m_j)$, on peut supposer $m_1\geq m_2\geq\dots
\geq m_j$, si bien que pour tout $i\leq j$, $m_i\leq j/i$.
On en d\'eduit la majoration
\[
\sup_{\substack{m\leq n \\ j\leq D}}\log^+ \abs{a_m^{(j)}}_v 
 \leq \sum_{j=1}^D \sup_{r\leq n/j} \log^+\abs{a_r}_v.
\]
d'o\`u finalement
\begin{equation} \label{eq.shidlovski}
 \frac 1n h_v(\phi^{(n)} \leq \sum_{j=1}^D \frac1n 
         \sup_{m\leq n/j} \log^+\abs{a_m}_v.
\end{equation}

Soit $T_1$ une partie finie de $T$ contenant les places archim\'ediennes.
En d\'ecomposant la somme sur les places $v\in T$ suivant $T_1$
et son compl\'ementaire, on obtient
\[ \limsup \frac1n \sum_{v\in T} h_v(\phi^{(n)})
  \leq \sum_{v\in T_1} \log^+ R_v(y)^{-1} 
         + \sum_{j=1}^D \frac1j
           \limsup \frac1n \sum_{v\in T\setminus T_1} \sup_{m\leq n}\log^+
\abs{a_m}_v, \]
c'est-\`a-dire
\[ \limsup \frac1n \sum_{v\in T} h_v(\phi^{(n)})
       \leq \rho_{T_1}(y)  +  \left(\sum_{j=1}^D\frac1j\right)
            \sigma_{T\setminus T_1} (y). \]
Prenant $T_1$ arbitrairement grand et utilisant le fait
que $\tau(y)=\inf_S \sigma_{S}(y)=0$, on obtient la majoration voulue.

\subsection{}
D{\'e}montrons maintenant le th{\'e}or{\`e}me d'alg{\'e}bricit{\'e}.
Soit $S\subset\Sigma_K$ 
un ensemble fini de places
contenant les places archim{\'e}diennes.
Par souci d'all{\`e}gement, on note
$E=E_{d,D}$. On note aussi $C'_v=\log \rg E$ pour $v$ archim{\'e}dienne
et $C'_v=0$ pour $v$ finie.

Compte tenu du lemme~\ref{Lemm.evidente}
et de la majoration~\eqref{eq.schwarz-simul}
pour toute place $v\in S$, 
l'in{\'e}galit{\'e} de pentes~\eqref{eq.inegpentes} entra\^{\i}ne
\begin{multline*}
 0 \leq \sum_{v\in S} \sum_{n=0}^\infty \rg(E^{(n)}/E^{(n+1)})
 \left( -n \log R'_v + (d+D)C_v +  C'_v\right) \\
  + C_S(\eps) \rg E + \sum_{n=0}^\infty  n\rg(E^{(n)}/E^{(n+1)})
 ( \rho_S(y)+ \eps(1+\log D) ) ,
\end{multline*}
d'o{\`u}, utilisant que $\sum_n \rg(E^{(n)}/E^{(n+1)}) =\rg E$,
\begin{multline} \label{ineq.pentes-bdw}
 \sum_{n=0}^\infty  n \rg(E^{(n)}/E^{(n+1)})
 \log \prod_{v\in S} R'_v
\\ 
\leq
\left( (d+D)\sum_{v\in S} C_v + C_S(\eps) + [K:\Q]\log\rg E \right)\rg E \\
 + \sum_{n=0}^\infty n\rg (E^{(n)}/E^{(n+1)})
        (\rho_S(y)+\eps(1+\log D)).
\end{multline}

Nous pouvons minorer $\Delta=\sum n \rg(E^{(n)}/E^{(n+1)})$
de la fa\c{c}on suivante.
Par construction, $\rg (E^{(n)}/E^{(n+1)})\leq \rg G^{(n)}=1$.
Il en r{\'e}sulte que $\rg E^{(n)}\geq \rg E -n$
et par suite,
\begin{align}
\sum_{n=0}^\infty n \rg(E^{(n)}/E^{(n+1)})
& = \sum_{n=1}^\infty \rg E^{(n)} \notag \\ 
& \geq \sum_{n=1}^{\rg E} (\rg E-n) = \frac12 \rg E (\rg E-1).
\end{align}
 
Faisons maintenant tendre $d$ vers $+\infty$, $D$ restant fixe.
On a ainsi $\Delta \gg d^2D^2$, si bien que 
$\rg E \log \rg E =\mathrm o(\Delta)$
et 
\[ \limsup \frac{(d+D) \rg E }{ \Delta} \leq \frac2{D+1}. \]
En divisant les deux membres de l'in\'egalit\'e~\eqref{ineq.pentes-bdw}
par $\sum n\rg(E^{(n)}/E^{(n+1)})$,
on obtient l'in\'egalit\'e
\[
\log \prod_{v\in S} R'_v \leq \rho_S(y)+\eps(1+\log D) + \frac1D \sum_{v\in
S} C_v. 
\]
On peut alors faire tendre $\eps$ vers~$0$, puis $D$ vers l'infini,
et enfin $R'_v$ vers $R_v$
dans cette in{\'e}galit{\'e} et l'on obtient
\begin{equation}
\log \prod_{v\in S} R_v \leq \rho_S(y). 
\end{equation}
Comme $\rho(y)<\infty$, 
le membre de droite peut {\^e}tre rendu
arbitrairement petit quitte {\`a} augmenter $S$,
mais ceci contredit alors l'hypoth{\`e}se que $\prod_{v\in\Sigma_K}R_v>1$.

\subsection{}
Pour d{\'e}montrer la seconde partie du th{\'e}or{\`e}me, on fixe encore
un ensemble fini $S\subset\Sigma_K$ contenant les places
archim{\'e}diennes.
On raisonne par l'absurde en supposant que $y$ n'est pas rationnelle.
Si $D=1$, l'application lin{\'e}aire $\phi$ est donc injective pour
tout entier $d\geq 1$.
On introduit un nouveau param{\`e}tre $N$ et on utilise la majoration
de $h_v(\phi^{(n)})$ fournie par 
le lemme~\ref{Lemm.evidente} si $v\not\in S$ ou si $n < N$
et par
l'estim{\'e}e~\eqref{eq.schwarz-triv} sinon.
On obtient ainsi l'in{\'e}galit{\'e}
\begin{align} \label{ineq.bd2}
0 & \leq C_{\complement S}(\eps) \rg E + (\rho_{\complement S}(y)+\eps)
            \sum_{n=0}^{\infty} n \rg(E^{(n)}/E^{(n+1)}) \notag\\
&\qquad {} +  \sum_{v\in S} \sum_{n\geq N}
         \left( (d-n)\log R'_v + D C_v + C'_v \right)
\rg(E^{(n)}/E^{(n+1)}) \\
&\qquad {} + \sum_{n < N} \rg(E^{(n)}/E^{(n+1)})
     (C_S(\eps)+[K:\Q]\log\rg E+2[K:\Q]\log n)  \notag\\
& \qquad{} 
       + (\rho_S(y)+\eps)\sum_{n< N} n\rg E^{(n)}/E{(n+1)} . \notag
\end{align}
Compte tenu de la majoration $\rg(E^{(n)}/E^{(n+1)})\leq 1$,
on a la majoration
\[
\Delta_N = \sum_{n < N} n\rg(E^{(n)}/E^{(n+1)})
 \leq \sum_{n < N} n = N(N-1)/2. \]
On fixe un param\`etre $\lambda$ et on pose $N=\lfloor 2\lambda d\rfloor$.
Lorsque $d$ tend vers l'infini,
mais bien s{\^u}r, $D=1$,
de sorte que $\rg E=2(d+1)$ et $N\sim \lambda\rg E$.
Puisque
\[ \Delta=\sum_{n=0}^\infty n\rg (E^{(n)}/E^{(n+1)})
         \geq \frac12 \rg E (\rg E-1), \]
on a $\limsup (\Delta_N/\Delta )\leq \lambda^2$.
De plus, si $\lambda\leq 1/2$, on a $N\leq d$
et on constate que $\rg E^{(N)}=\rg E-N$, si bien que
\[
 \limsup_{d\ra +\infty} \frac{d\rg E^{(N)}}{\Delta} 
       \leq 1-\lambda.
\]
Apr\`es division des deux membres par $\Delta$ et passage
\`a la limite, l'in\'egalit\'e~\eqref{ineq.bd2} devient
\[ 0 \leq (\rho_{\complement S}(y)+\eps) + \lambda^2 (\rho_S(y)+\eps)
        + (\lambda^2-\lambda) \log\prod_{v\in S} R'_v , \]
d'o\`u, si $\eps$ tend vers~$0$ et $R'_v$ vers $R_v$,
\[ \lambda(1-\lambda) \log\prod_{v\in S}R_v
 \leq \rho_{\complement S}(y) + \lambda^2 \rho_S(y). \]
Comme $\rho_S(y)\leq \rho(y)$,
cette in\'egalit\'e contredit l'hypoth\`ese $\prod_{v\in \Sigma_K} R_v>1$,
lorsque  $S$ est assez grand.

(Il n'est pas certain que cette 
d{\'e}monstration du th{\'e}or{\`e}me de Borel-Dwork
soit plus simple que la d{\'e}monstration originelle, cf.~\cite{eborel1894},
\cite{dwork60} ou~\cite{amice75}.)

\section{Alg{\'e}bricit{\'e} de sous-vari{\'e}t{\'e}s formelles}\label{sec.bost}

\subsection{}
Nous voulons maintenant indiquer
la d{\'e}monstration du th{\'e}or{\`e}me~\ref{theo.feuilletage}.
Rappelons la situation: $X$ est une vari{\'e}t{\'e} alg{\'e}brique
lisse sur un corps de nombres $K$, $P$ est un point
de $X(K)$  et $F$ est un sous-fibr{\'e} involutif de $\mathrm TX$
qui, modulo~$p$ pour presque tout nombre premier $p$, est stable
par puissance~\pieme.
On suppose de plus qu'il existe un plongement
$\sigma_0\colon K\hra\C$ tel que la feuille passant par $P$ du feuilletage
holomorphe de $X_{\sigma_0}(\C)$ induit par $F$ est {\og uniformis{\'e}e par 
l'espace affine\fg}. C'est une hypoth{\`e}se plus contraignante que
celle du th{\'e}or{\`e}me~\ref{theo.feuilletage}
mais suffisante pour {\'e}tablir le th{\'e}or{\`e}me~\ref{theo.lie}.

On veut alors en d{\'e}duire que ladite feuille est une sous-vari{\'e}t{\'e}
alg{\'e}brique de $X$, ou encore, si $\widehat V$ d{\'e}signe la feuille
formelle de $F$ en $P$, qu'il existe une sous-vari{\'e}t{\'e} alg{\'e}brique
$V$ de $X$ dont $\widehat V$ est le compl{\'e}t{\'e} en $P$.

Soit $Y$ l'adh{\'e}rence de Zariski de $\widehat V$ dans $X$,
c'est-{\`a}-dire la plus petite sous-vari{\'e}t{\'e} alg{\'e}brique de $X$
dont le compl{\'e}t{\'e} en $P$ contient $\widehat V$.
Il suffit de d{\'e}montrer que $\dim Y=\dim\widehat V$
car cette hypoth{\`e}se implique que $Y$ est automatiquement 
lisse et une feuille du flot~$F$.
On peut aussi supposer que $X$ est projective (mais $X$ n'est
alors lisse que dans un voisinage de $P$); $Y$ est alors projective.
Soit $L$ un fibr{\'e} inversible ample sur $X$.
L'hypoth{\`e}se que $Y$ est l'adh{\'e}rence de $\widehat V$ signifie
que pour tout entier $D\geq 0$,
l'homomorphisme de restriction {\`a} $\widehat V$,
\begin{equation}
\phi_D \colon \Gamma(Y,L^D) \ra \Gamma(\widehat V, L^D)
\end{equation}
est injectif.

Sous l'hypoth{\`e}se que $\widehat V$ n'est pas alg{\'e}brique, c'est-{\`a}-dire
que $\dim Y>\dim\widehat V$, nous voulons contredire cette injectivit{\'e}.
Nous allons pour cela utiliser l'in{\'e}galit{\'e} de pentes 
(prop.~\ref{prop.pentes}).
Pla\c{c}ons-nous ainsi dans le contexte du chapitre~\ref{sec.arakelov}.
Choisissons des mod{\`e}les entiers de toute la situation:
\begin{itemize}
\item un $\mathfrak o_K$-sch{\'e}ma propre et plat~$\mathscr X$ 
tel que $\mathscr X\otimes K=X$ ;
\item 
une section $\eps_P\colon \Spec\mathfrak o_K\ra\mathscr X$
prolongeant $P\in X(K)$ ;
\item  l'adh{\'e}rence sch{\'e}matique $\mathscr Y$
de $\widehat V$ (ou de $Y$) dans $\mathscr X$ ;
\item 
un fibr{\'e} inversible $\mathscr L$ sur $\mathscr X$
dont la restriction {\`a} $X$ est {\'e}gale {\`a} $L$.
\end{itemize}
On note $\mathrm t^\vee$ l'image de $\eps_P^*\Omega^1_{\mathscr Y/\mathfrak
o_K}$ dans $\Omega^1_{Y/K,P}$ ; c'est un $\mathfrak o_K$-module
projectif de rang~$\dim Y$.
Choisissons aussi 
\begin{itemize}
\item une m{\'e}trique hermitienne sur le fibr{\'e} holomorphe
induit par $L$ sur $X(\C)$ ;
\item une mesure de Lebesgue positive sur $Y(\C)$ ;
\item une m{\'e}trique hermitienne sur l'espace tangent en $P$
{\`a} $\widehat V$, ou par dualit{\'e}, sur $\mathrm t^\vee$ ;
\end{itemize}
toutes invariantes par la conjugaison complexe.

\subsection{}
Pour tout entier $D\geq 1$,
on notera $\bar E_D=\Gamma(\mathscr Y,L^D)$, muni de sa structure naturelle
de $\mathfrak o_K$-fibr{\'e} vectoriel hermitien d{\'e}finie comme au
paragraphe~\ref{arak.globalsect}.
Le $K$-espace vectoriel $\Gamma(\widehat V,L^D)$ est 
filtr{\'e} par l'ordre d'annulation en~$P$, les sous-quotients
successifs s'identifient aux fibres g{\'e}n{\'e}riques des $\mathfrak o_K$-fibr{\'e}s
vectoriels hermitiens
\begin{equation}
      \Sym^n \bar {\mathrm t}^\vee \otimes_{\mathfrak o_K} \eps_P^*\bar{\mathscr
L}^D.
\end{equation}
Si $\bar E^{(n)}_D$ d{\'e}signe l'image inverse de cette filtration
par l'homomorphisme d'{\'e}valuation $\phi_D$, l'application
\begin{equation}
  \phi_D^{(n)}\colon   E^{(n)}_D/E^{(n+1)}_D \ra  \Sym^n {\mathrm t}_K^\vee
\otimes L_P^D
\end{equation}
s'identifie {\`a} l'application {\og jet d'ordre~$n$ en $P$\fg}.
L'homomorphisme $\bar E_D\ra \Gamma(\widehat V,L^D)$
est injectif par construction.
L'in{\'e}galit{\'e} de pentes de la proposition~\ref{prop.pentes} s'{\'e}crit alors
\begin {equation}\label{ineq.pentes}
\hdeg \bar E_D
\leq \sum_{n\geq 0} \rg(E_D^{(n)}/E_D^{(n+1)})
 \left( \hmumax( \Sym^n \bar {\mathrm t}^\vee \otimes_{\mathfrak o_K}
\eps_P^*\bar{\mathscr L}^D ) 
      + h(\phi_D^{(n)}) \right)
.
\end{equation}

Dans~\cite{bost2001}, $h(\phi_D^{(n)})$ est major{\'e}e par la proposition
suivante.
\begin{Prop}\label{prop.unif1}
Rappelons l'hypoth{\`e}se:
\begin{itemize}
\item pour presque tout nombre premier $p$,
la r{\'e}duction de $F$ modulo~$p$ est stable par puissance~\pieme;
\item il existe une vari{\'e}t{\'e} complexe $M$ v{\'e}rifiant la propri{\'e}t{\'e}
de Liouville, un point $O\in M$,
un plongement complexe $\sigma_0\colon K\hra\C$
et une application holomorphe $\psi$ de $M$ vers la feuille en $P$
du feuilletage holomorphe induit par $F$ sur $X_{\sigma_0}(\C)$
telle que $\psi(O)=P$ et telle que $\psi$ soit biholomorphe
d'un voisinage de $O$ vers un voisinage de $P$ dans cette feuille.
\end{itemize}
Alors, pour tout $\rho>0$, il existe un r{\'e}el $C(\rho)$ tel
que pour tous les entiers $D$ et $n\geq 1$, on ait
\begin{equation}\label{ineq.liouville}
h(\phi_D^{(n)}) \leq -n \rho + D C(\rho) .
\end{equation}
% Alors, il existe une fonction $\tau\colon \R_+\ra\R$ minor{\'e}e
% telle que  
% \begin{equation}
% \lim_{x\ra +\infty} \frac{\tau(x)}{x} = +\infty
% \end{equation}
% et telle que pour tous $n$ et $D$, on ait 
% \begin{equation}
% h(\phi_D^{(n)}) \leq -D \tau (n/D).
% \end{equation}
\end{Prop}

Pour d{\'e}montrer le th{\'e}or{\`e}me~\ref{theo.lie},
il suffit de traiter le cas o{\`u} $M=\C^d$ dans lequel l'analyse
aux places archim{\'e}diennes peut {\^e}tre pr{\'e}sent{\'e}e de fa\c{c}on
relativement {\'e}l{\'e}mentaire (voir le
paragraphe~\ref{subsec.liouville}).

\subsection{}
Tirons maintenant la contradiction de ces estimations.
En ins{\'e}rant
dans l'in{\'e}galit{\'e}~\eqref{ineq.pentes} 
la majoration~\eqref{ineq.liouville}, les
conclusions de la proposition~\ref{prop.hsa} et du lemme~\ref{lemm.symk}
ainsi que la formule~\eqref{eq.hmumax(etensl)},
on obtient l'in{\'e}galit{\'e}, valable pour tout $\rho>0$,
\begin{equation}
 - c_1 D \rg E_D \leq \sum_{n\geq 0}\rg(E_D^{(n)}/E_D^{(n+1)})
       \left( nc_2+ D c_3 + D C(\rho)-n\rho \right) ,
\end{equation}
soit encore
\begin{equation}\label{ineq.pente3}
(\rho-c_2) \sum n\rg (E_D^{(n)}/E_D^{(n+1)})
     \leq  (c_3+c_1 + C(\rho)) D \rg E_D.
\end{equation}
Or, d'apr{\`e}s l'alg{\`e}bre lin{\'e}aire, si $d=\dim \widehat V$,
\begin{gather}
\rg (E_D^{(n)}/E_D^{(n+1)}) 
         \leq \rg \Sym^n \mathrm t^\vee = \binom{n+d-1}{d-1} \\
\rg E_D^{(n)} \geq \rg E_D - \sum_{k=0}^{n-1} \binom{k+d-1}{d-1}
        = \rg E_D - \binom{n+d-1}{d} 
\end{gather}
tandis que le th{\'e}or{\`e}me de Hilbert garantit que
\begin{equation}
\rg E_D \simeq c_4 D^{\dim Y}.
\end{equation}
Il en r{\'e}sulte que, si $N$ est un entier $\geq 1$ arbitraire,
\begin{align*}
\sum_{n=0}^\infty n \rg (E_D^{(n)}/E_D^{(n+1)}) & =
 \sum_{n=1}^\infty \rg E_D^{(n)} \geq \sum_{n=1}^N \rg E_D^{(n)} \\
&\geq N \rg E_D -  \binom{N+d}{d+1} .
% \geq c_4 N D^{\dim Y} - c_5 N^{d+1}.
\end{align*}
Supposons par l'absurde que $\dim Y>d$ et
choisissons $N=\lfloor D^{\alpha}\rfloor$ pour un r{\'e}el
$\alpha$ tel que $1<\alpha<\dim Y/d$. Il en r{\'e}sulte
que $D=\mathrm o(N)$ et $N^d=\mathrm o(D^{\dim Y})$. Par suite,
$ \binom{N+d}{d+1}= \mathrm o(N \rg E_D)$.
En passant {\`a} la limite dans l'in{\'e}galit{\'e}~\eqref{ineq.pente3},
on obtient alors
\begin{equation}
\rho -c_2 \leq 0.
\end{equation}
Comme $\rho$ est arbitraire, on a une contradiction.

\subsection{}\label{subsec.liouville}
Il reste {\`a} {\'e}tablir l'in{\'e}galit{\'e}~\eqref{ineq.liouville},
et pour cela, nous allons majorer $h_v(\phi_D^{(n)})$ pour toute
place $v$ de $\Sigma_K$.
Lorsque $v$ est une place finie, il s'agit essentiellement
de {\og contr{\^o}ler les d{\'e}nominateurs\fg} 
qui apparaissent dans le th{\'e}or{\`e}me de Frobenius.

Notons $s$ la dimension de~$X$ et
fixons des coordonn{\'e}es locales (pour la topologie {\'e}tale)
$x_1,\dots,x_s$ sur $X$ au point~$P$, de sorte qu'autour
de $P$, $\mathrm TX$ est libre de base $\frac{\partial}{\partial x_1}$,
\dots, $\frac{\partial}{\partial x_s}$. Les $x_i$
d{\'e}finissent aussi un isomorphisme formel
$\widehat\alpha_0 \colon \widehat X_P \simeq \widehat\A^s$.
Quitte {\`a} renum{\'e}roter les $x_i$, il est possible de trouver 
une base $D_1,\dots,D_d$ des sections de $F$ dans un voisinage de $P$ de
la forme
\begin{equation}
D_i = \frac{\partial}{\partial x_i} + \sum_{j=d+1}^s a_{i,j}
\frac{\partial}{\partial x_j} ,
\end{equation}
o{\`u} les $a_{i,j}$ sont des s{\'e}ries formelles.
Soit $\widehat\alpha\colon \widehat X_P \ra \widehat\A^d$ le morphisme
d{\'e}duit de $\widehat\alpha_0$ par projection sur les $d$ premi{\`e}res coordonn{\'e}es.
On a ainsi $\widehat\alpha_*(D_i)=\frac{\partial}{\partial x_i}$ et
\[ \widehat\alpha_*([D_i,D_j]) = [\widehat\alpha_*(D_i),\widehat\alpha_*(D_j)] = [ \frac{\partial}{\partial x_i},
 \frac{\partial}{\partial x_j} ] = 0 \]
si bien que $[D_i,D_j]$ est une combinaison lin{\'e}aire
des $\frac{\partial}{\partial x_k}$ pour $k>d$.
Comme $F$ est involutif, $[D_i,D_j]\in F$.
Il en r{\'e}sulte que $[D_i,D_j]=0$: les champs de vecteurs
$D_1,\dots,D_d$ commutent.

Ainsi, $\widehat V$ est la feuille passant par~$P=(0,\dots,0)$
du {\og flot formel\fg},
$\Phi\colon \widehat {\A}\vphantom{\A}^d \times \widehat{\A}\vphantom{\A}^{s}
\ra  \widehat{\A}\vphantom{\A} ^{s} $,
\begin{multline}
   ((t_1,\dots,t_d),(x_1,\dots,x_s)) \mapsto \exp\big(\sum_{i=1}^d t_i
e_i\big)
 \cdot (x_1,\dots,x_s)  \\
         =\sum_{n_1,\dots,n_d\geq 0} \frac{{t_1}^{n_1}}{n_1!}
 \cdots \frac{t_d^{n_d}}{n_d!}  {e_1}^{n_1}\cdots e_d^{n_d}
\cdot(x_1,\dots,x_s).
\end{multline}
Pour tout multi-indice $n=(n_1;\dots;n_d)$ et
tout entier $i\in\{1;\dots;s\}$,  notons
$  P_{n,i}$ la s{\'e}rie formelle
$  \frac{1}{n_1!\dots n_d!} e_1^{n_1}\cdots e_d^{n_d}  (x_i)$,
de sorte que $\widehat V$ admet la param{\'e}trisation formelle
\[ (t_1,\dots,t_d) \mapsto \sum_{n\in\mathbf N^d} 
  t_1^{n_1}\cdots t_d^{n_d} (P_{n,1}(0),\dots,P_{n,s}(0)) .
\]
% Ce flot d{\'e}finit un isomorphisme
% $\phi\colon  \widehat\A\vphantom{\A}^s \ra
% \widehat V \times \widehat \A\vphantom{\A}^{s-d}$,
% \begin{equation}
%  (t_1,\dots,t_s) \ra \psi( (t_1,\dots,t_d),(0,\dots,0))
%  + (0,\dots,0,t_{d+1},\dots,t_s)
% \end{equation}
% dont la diff{\'e}rentielle en l'origine est
% \begin{equation}
% \mathrm d\phi = e_1\,\mathrm dt_1+\dots+e_d\,\mathrm dt_d
%      + \mathrm dt_{d+1}+\dots+\mathrm dt_s.
% \end{equation}

\begin{Lemm}
Pour toute place finie $v$ de $K$, de caract{\'e}ristique
r{\'e}siduelle $p$, il existe un r{\'e}el $C_v\geq 0$
tel que  pour tous $n\in\mathbf N^d$ et tout $i\in\{1;\dots;s\}$,
\[ - \log \abs{P_{n,i}(0)}_v \leq (n_1+\dots+n_d) C_v. \]
De plus, pour presque tout $v$, on peut choisir
\[ C_v = [K:\Q] \frac{\log p}{p(p-1)}. \]
\end{Lemm}
L'existence d'un tel $C_v$ {\'e}quivaut {\`a} la convergence $v$-adique
du flot formel dans un polydisque.

De plus, pour presque toute place finie $v$,
les coordonn{\'e}es
locales $x_i$ s'{\'e}tendent en des coordonn{\'e}es locales
sur $\mathscr X \otimes \mathfrak o_v$, ainsi que les champs
de vecteurs locaux $D_1,\dots,D_d$.
Si de plus la r{\'e}duction modulo l'id{\'e}al premier $\mathfrak p_v$ 
est stable par puissance~\pieme,
on constate que pour tout $i$, en r{\'e}duction modulo $\mathfrak p_v$,
\[ \widehat\alpha_*(D_i^p)  =\widehat\alpha_*(D_i)^p= 
\frac{\partial^p}{\partial x_i^p} = 0. \]
Comme la r{\'e}duction modulo $\mathfrak p_v$ de $F$
est suppos{\'e}e stable par puissance~\pieme,
il en r{\'e}sulte que $D_i^p=0$ modulo~$\mathfrak p_v$.
Si de plus $p$ ne divise pas le discriminant de $K$,
on a alors $v(p)=1$ et
\[ 
v( e_1^{n_1}\dots e_d^{n_d}\cdot x_i )
       \geq \lfloor  n_1/p\rfloor + \dots + \lfloor n_d/p\rfloor
\]
et 
\[ v(P_{n,i}(0)) \geq 
       - \sum_{j=1}^d \left(v(n_j!) - \lfloor n_j/p\rfloor \right)
\geq - \frac{n_1+\dots+n_d}{p(p-1)} \]
d'o{\`u} l'on d{\'e}duit le lemme.

De ces estim{\'e}es
d{\'e}coule facilement (lemme de Schwarz $v$-adique) une majoration,
valable pour toute place $v$ finie,
\begin{gather} \label{ineq.hv}
     h_v(\phi_D^{(n)}) \leq n C_v \\
\intertext{et de plus, }
\label{ineq.hv2}
\sum_{v\in\Sigma_{K,f}} C_v<+\infty.
\end{gather}

Pour tout plongement complexe $\sigma\colon K\hra\C$,
soit $M_\sigma=B(0,R_\sigma)^d$ le polydisque ouvert de centre~$O$ et
de rayon $R_\sigma\in\R_+^*\cup\{+\infty\}$
et soit $\psi_\sigma\colon M_\sigma\ra \mathscr Y_\sigma(\C)$
une application holomorphe induisant un isomorphisme
d'un voisinage de $O$ dans $M_\sigma$ avec un voisinage
de $P$ dans la feuille holomorphe passant par $P$ du feuilletage
d{\'e}fini par $F$.
Par hypoth{\`e}se, de telles applications existent pour tout 
$\sigma$, et pour la place $\sigma_0$, on a $R_{\sigma_0}=+\infty$.

Fixons une section globale sans z{\'e}ro $\eps_\sigma$ de  $\psi_\sigma ^*L$.
Si $s\in E_D^n$, on {\'e}crit $\psi_\sigma^* s = f \eps_\sigma^D$
o{\`u}~$f$ est une fonction holomorphe sur $M_\sigma$
s'annulant {\`a} l'ordre~$n$ en l'origine.

Soit pour tout $\sigma$ un r{\'e}el $R'_\sigma$ tel que
$0<R'_\sigma<R_\sigma$.
Soit $j^nf$ le jet d'ordre~$n$ de $f$ en l'origine:
c'est l'application polynomiale homog{\`e}ne de degr{\'e}~$n$:
\[ \C^d \ra \C, \qquad (a_1,\dots,a_d) \mapsto \sum_{n_1+\dots+n_d=n}
          \frac{ a_1^{n_1}}{n_1!} \dots \frac{a_d^{n_d}}{n_d!} 
                 \frac{\partial ^{n_1}}{\partial z_1^{n_1}}\dots
\frac{\partial ^{n_d}}{\partial z_d^{n_d}} f (0). \]
D'apr{\`e}s le lemme de Schwarz usuel,
on a pour tout $(a_1,\dots,a_d)\in\C^d$ l'in{\'e}galit{\'e}
\[ \abs{ j^nf(a_1,\dots,a_d) } \leq \max(\abs{a_i})^n \norm{f}_{R'_\sigma}
 (R'_\sigma)^{-n},  \]
$\norm{f}_{R'_\sigma}$ d{\'e}signant le sup de $\abs{f}$
sur le polydisque $\bar B(0,R'_\sigma)^d$
et  la norme du jet de $s$ est major{\'e}e par
\[ \norm{j^n s} \leq \norm{\mathrm d\psi_{\sigma}^{-1}}^{n}
      (R'_\sigma)^{-n} \norm{\psi^* s}_{R'_\sigma}
        \norm{\eps_\sigma^{-1}}_{R'_\sigma}^{D}. \]
Finalement, compte tenu de l'in{\'e}galit{\'e}~\eqref{lemm.gromov},
il existe deux constantes $A$ et $B(R'_\sigma)$  telles que
\[ \norm{j^n s} \leq (R'_\sigma/A)^{-n} B(R'_\sigma)^D \]
et
\begin{equation}\label{ineq.hsigma}
 h_\sigma(\phi_D^{(n)}) \leq D \log B(R'_\sigma) -n \log (R'_\sigma/A).
\end{equation}
L'in{\'e}galit{\'e}~\eqref{ineq.liouville}
d{\'e}coule alors imm{\'e}diatement des in{\'e}galit{\'e}s~\eqref{ineq.hv}, \eqref{ineq.hv2}
et~\eqref{ineq.hsigma} et du fait que $R'_{\sigma_0}$
puisse {\^e}tre pris arbitrairement grand.

\begin{Rema}
La majoration~\eqref{ineq.hv2} est un analogue de la condition
$\tau=0$ du chapitre~\ref{sec.abdw}. Ce n'est d'ailleurs
pas surprenant: dans le contexte de la conjecture de Grothendieck,
la condition $\tau=0$ provient justement 
de l'hypoth{\`e}se d'annulation des $p$-courbures, de m{\^e}me
que l'in{\'e}galit{\'e}~\eqref{ineq.hv2} a {\'e}t{\'e} {\'e}tablie en utilisant
la $p$-int{\'e}grabilit{\'e} du feuilletage.
\end{Rema}

\backmatter

\providecommand{\noopsort}[1]{}\providecommand{\url}[1]{\textit{#1}}
\providecommand{\bysame}{\leavevmode ---\ }
\providecommand{\og}{``}
\providecommand{\fg}{''}
\providecommand{\smfandname}{\&}
\providecommand{\smfedsname}{\'eds.}
\providecommand{\smfedname}{\'ed.}
\providecommand{\smfmastersthesisname}{M\'emoire}
\providecommand{\smfphdthesisname}{Th\`ese}

\end{document}